# EFFICIENT ESTIMATION OF BANACH PARAMETERS IN SEMIPARAMETRIC MODELS[1]

By Chris A. J. Klaassen and Hein Putter

*University of Amsterdam and University of Leiden*

Consider a semiparametric model with a Euclidean parameter and an infinite-dimensional parameter, to be called a Banach parameter. Assume:

(a) There exists an efficient estimator of the Euclidean parameter.

(b) When the value of the Euclidean parameter is known, there exists an estimator of the Banach parameter, which depends on this value and is efficient within this restricted model.

Substituting the efficient estimator of the Euclidean parameter for the value of this parameter in the estimator of the Banach parameter, one obtains an efficient estimator of the Banach parameter for the full semiparametric model with the Euclidean parameter unknown. This hereditary property of efficiency completes estimation in semiparametric models in which the Euclidean parameter has been estimated efficiently. Typically, estimation of both the Euclidean and the Banach parameter is necessary in order to describe the random phenomenon under study to a sufficient extent. Since efficient estimators are asymptotically linear, the above substitution method is a particular case of substituting asymptotically linear estimators of a Euclidean parameter into estimators that are asymptotically linear themselves and that depend on this Euclidean parameter. This more general substitution case is studied for its own sake as well, and a hereditary property for asymptotic linearity is proved.

**1. Introduction.** Estimation of a parameter is not a goal in itself. Typically, the purpose is to determine a reliable picture of future behavior of a random system. In semiparametric models this means that estimation of just the finite-dimensional, Euclidean parameters does not finish the job.

Received August 1999; revised April 2004.

[1]Supported in part by the Netherlands Organization for the Advancement of Scientific Research (NWO), via the project Computationally Intensive Methods in Stochastics.

The values of the Banach parameters are needed to complete the picture. The situation in classical parametric models is similar. Consider linear regression under normal errors with unknown variance. The regression parameters are the parameters of interest, but the variance of the errors, although of secondary interest, is essential to describe the behavior of the dependent variable at a particular value of the independent variable, as for instance in prediction. In semiparametric linear regression the error distribution with mean zero is completely unknown. Again, this distribution is essential in describing the behavior of the dependent variable. Therefore, its distribution function has to be estimated. This may be done along the following lines:

1. Estimate the regression parameter vector $\theta$ efficiently using (by now standard) semiparametric theory.
2. Given the true value of the parameter $\theta$, the error distribution function $G$ can be estimated efficiently, since the i.i.d. errors can be reconstructed from the observations in this case.
3. Using the estimated value of $\theta$, construct the residuals and instead of the i.i.d. errors use these residuals to estimate the Banach parameter $G$ in the same way as in step 2.

The crux of the present paper is that the resulting estimator of $G$ is efficient. In fact, for any semiparametric model, we will prove that this approach, which is in line with statistical practice, yields an efficient estimator of the Banach parameter, provided a sample splitting scheme is applied. Since we assume that efficient estimators of $\theta$ are available, we shall focus on efficient estimation of the Banach parameter $G$ in the presence of the Euclidean nuisance parameter $\theta$. Sample splitting is unnecessary and the direct substitution estimator works if the conditional estimator of the Banach parameter given $\theta$ depends on $\theta$ in a smooth way. In order to be able to estimate $G$ efficiently according to our approach, it is essential in non-adaptive cases that in step 1 the Euclidean parameter $\theta$ can be estimated efficiently in the semiparametric sense. The Banach parameter needed for more complete inference, like the distribution function $G$ of the errors in semiparametric linear regression, typically is *unequal* to the Banach parameter needed in efficient semiparametric estimation of $\theta$, this parameter being the score function $-d\log(dG(x)/dx)/dx$ for location in the linear regression model. In fact, Klaassen (1987) has shown that $\theta$ can be estimated efficiently if and only if the efficient influence function for estimating $\theta$ can be estimated consistently and $\sqrt{n}$-unbiasedly, given $\theta$, and $\theta$ can be estimated $\sqrt{n}$-consistently, with $n$ denoting sample size. Of course, this efficient influence function depends on the Banach parameter of interest, but typically differs from it.

To give a more explicit and precise statement of our results, let $\mathcal{P}$ be our semiparametric model given by

$$(1.1) \qquad \mathcal{P} = \{P_{\theta,G} : \theta \in \Theta, G \in \mathcal{G}\}, \qquad \Theta \subset \mathbb{R}^k, \mathcal{G} \subset \mathcal{H},$$

where $\Theta$ is an open subset of $\mathbb{R}^k$ and $\mathcal{G}$ is a subset of a Banach or preferably a Hilbert space $(\mathcal{H}, \langle \cdot, \cdot \rangle_{\mathcal{H}})$. Typically, in a natural parametrization, $G$ would be a distribution function and hence an element of a Banach space $L_\infty$. If a $\sigma$-finite measure would dominate the distributions in $\mathcal{G}$, then an obvious parametrization would be via the corresponding densities $g$ of $G$, which are elements of a Banach space $L_1$. However, via the square roots $\sqrt{g}$ we parametrize by elements of a Hilbert space $L_2$. Therefore, we shall assume that $\mathcal{G}$ is a subset of a Hilbert space or can be identified with it. We are interested in estimating a parameter

$$(1.2) \qquad \nu = \nu(P_{\theta,G}) = \tilde{\nu}(G),$$

where $\tilde{\nu} : \mathcal{G} \to \mathcal{B}$ ($\mathcal{B}$ Banach space) is pathwise differentiable; see (4.4) for details. Estimation has to be based on i.i.d. random variables $X_1, X_2, \ldots, X_n$ with unknown distribution $P_{\theta,G} \in \mathcal{P}$ on the sample space $(\mathcal{X}, \mathcal{A})$. Let $\mathcal{P}_2(\theta)$ be the submodel of $\mathcal{P}$ where $\theta$ is known. Let the submodel estimator $\hat{\nu}_{\theta,n}$ be an efficient estimator of $\nu$ within $\mathcal{P}_2(\theta)$. Suppose that we also have an estimator $\hat{\theta}_n$ of $\theta$ at our disposal within $\mathcal{P}$. Following step 3 above, an obvious candidate for estimating $\nu$ in the full model $\mathcal{P}$ would be the substitution estimator $\hat{\nu}_{\hat{\theta}_n,n}$. We shall show that a split-sample modification of $\hat{\nu}_{\hat{\theta}_n,n}$ is an efficient estimator of $\nu$ in $\mathcal{P}$ if $\hat{\theta}_n$ is an efficient estimator of $\theta$ in $\mathcal{P}$. In adaptive cases, for $\hat{\nu}_{\hat{\theta}_n,n}$ to be efficient in $\mathcal{P}$ it is sufficient that the estimator $\hat{\theta}_n$ be $\sqrt{n}$-consistent. The substitution estimator $\hat{\nu}_{\hat{\theta}_n,n}$ itself is semiparametrically efficient if the submodel estimator $\hat{\nu}_{\theta,n}$ depends smoothly on $\theta$, which is typically the case.

The asymptotic linearity of the efficient estimators involved warrants the resulting substitution estimator to be asymptotically linear as well. We study this hereditary property of asymptotic linearity of estimators for its own sake in Section 2, where we refrain from the efficiency assumptions made above. In Section 3 we discuss such simple examples as the sample variance and estimators of the standardized error distribution in linear regression. There we will also introduce models that we propose to call parametrized linkage models. In Section 4 we will collect some results about efficient influence functions in the various (sub)models that we consider. Section 5 contains our main results for efficiency. In Section 6 we will discuss a number of examples.

A general class of semiparametric models $\mathcal{P} = \{P_{\theta,G} : \theta \in \Theta, G \in \mathcal{G}\}$ in which our results apply is the class of models that can be handled by profile likelihood. If $l_n(\theta, G)$ is the appropriately defined likelihood of $n$ independent observations from $\mathcal{P}$, then a maximum likelihood estimator $\hat{\theta}_n$ of $\theta$ can be found by maximizing the *profile likelihood*

$$pl_n(\theta) = \sup_{G \in \mathcal{G}} l_n(\theta, G). \tag{1.3}$$

This amounts to maximizing the likelihood in two steps. First maximize with respect to $G$ for a given $\theta$. The maximizer of $l_n(\theta, G)$ with respect to $G$, say $\hat{G}_n(\theta)$, will generally depend on $\theta$. Placing the submodel estimator $\hat{G}_n(\theta)$ back into the likelihood, we obtain a function of $\theta$ only, $pl_n(\theta) = l_n(\theta, \hat{G}_n(\theta))$. Murphy and van der Vaart (2000) show that the profile likelihood can to a large extent be viewed as an ordinary likelihood. In particular, under some regularity conditions, asymptotic efficiency of the maximizer $\hat{\theta}_n$ of (1.3) can be proved. Important in this construction is the fact that the maximizer of the likelihood with respect to $G$, obtained in the first maximization step, $\hat{G}_n(\theta)$, is not yet a complete estimator of $G$. This submodel estimator is only an estimator of $G$ for a given value of $\theta$, just as in step 2 of our linear regression example above. Having found an efficient estimator of $\theta$, estimation of $(\theta, G)$ is then completed by considering the obvious substitution estimator $\hat{G}_n = \hat{G}_n(\hat{\theta}_n)$. The estimator $\hat{G}_n(\theta)$ for given $\theta$ is already available as a result of the maximizing step in (1.3). The Banach parameter $\tilde{\nu}(G)$ or $G$ itself will not generally be estimable at $\sqrt{n}$-rate, but it may be possible to estimate real-valued functionals $\kappa = \kappa(P_{\theta,G}) = \tilde{\kappa}(G)$ of $G$ at $\sqrt{n}$-rate. In cases where $\tilde{\kappa}(\hat{G}_n(\theta))$ is an efficient estimator of $\tilde{\kappa}(G)$ given $\theta$, our results can be applied to yield a fully efficient estimator $\tilde{\kappa}(\hat{G}_n(\hat{\theta}_n))$ of $\tilde{\kappa}(G)$. Numerous examples fall into this class, some of them treated in some detail in this paper, like the Cox proportional hazards model for right censored data (Example 6.6) and for current status data [Huang (1996) and Bolthausen, Perkins and van der Vaart (2002)], frailty models [Nielsen, Gill, Andersen and Sørensen (1992)], the proportional odds model [Murphy, Rossini and van der Vaart (1997)], selection bias models [Gilbert, Lele and Vardi (1999) and Cosslett (1981)] and random effects models [Butler and Louis (1992)].

We will consider a number of examples in more detail in Section 6, namely estimation of the variance with unknown mean, estimation of the error distribution in parametrized linkage models and in particular in the location problem with the bootstrap as an application, estimation of a (symmetric) error distribution in linear regression as an example of the adaptive case, and finally, estimation of the baseline distribution function in the Cox proportional hazards model.

The framework of the present paper has been presented in Klaassen and Putter (1997) within the linear regression model with symmetric error distribution and has been used by Müller, Schick and Wefelmeyer (2001) in their discussion of substitution estimators in semiparametric stochastic process models.

There are fundamental theorems of algebra, arithmetic and calculus. Statistics has its fundamental rule of thumb. It states that "replacing unknown parameters in statistical procedures by estimators of them yields appropriate procedures." This paper describes a large class of estimation problems where this rule of thumb is indeed a *theorem*.

**2. Heredity of asymptotic linearity of substitution estimators.** In this section we will study the local asymptotic behavior of estimators that are obtained by combining two asymptotically linear estimators in the way described in Section 1. We will prove the hereditary property that under certain regularity conditions the resulting estimators are asymptotically linear as well and we will describe their influence functions. The main application of this heredity result is to efficient estimators as described in Section 1. This will be pursued in Section 5, but we believe the hereditary property is of independent interest as well.

Although we will apply this hereditary property to semiparametric models $\mathcal{P}$ as in (1.1), we will be able to restrict attention in the present section to parametric models since the phenomenon under study occurs within the natural parametric submodels $\mathcal{P}_1(G) = \{P_{\theta,G} : \theta \in \Theta\}$ of $\mathcal{P}$ with $G \in \mathcal{G}$ fixed.

So, within this section, let $\mathcal{P} = \{P_\theta : \theta \in \Theta\}$, $\Theta \subset \mathbb{R}^k$ open, be a parametric model, and let $X_1, X_2, \ldots, X_n$ be the i.i.d. random variables with distribution $P_\theta \in \mathcal{P}$ on the sample space $(\mathcal{X}, \mathcal{A})$ that are used for estimation. Since our considerations are of the usual local asymptotic type, we introduce an arbitrary fixed $\theta_0 \in \Theta$ at which the local asymptotics is focused.

For every $m \in \mathbb{N}$ let $\Psi_m$ be the set of all measurable functions $\psi$ from $\mathcal{X} \times \Theta$ into $\mathbb{R}^m$ such that $\int \psi(x;\theta) \, dP_\theta(x) = 0$ and $\int |\psi(x;\theta)|^2 \, dP_\theta(x) < \infty$ for all $\theta \in \Theta$, where $|\cdot|$ denotes a Euclidean norm. Fix $m \in \mathbb{N}$ and consider a differentiable function $\kappa$ from $\Theta$ into $\mathbb{R}^m$.

DEFINITION 2.1. An estimator $\hat{\kappa}_n$ of $\kappa(\theta)$ is *locally asymptotically linear* at $\theta_0$ if there exists a $\psi \in \Psi_m$ such that

$$\sqrt{n}\left|\hat{\kappa}_n - \kappa(\theta_n) - n^{-1}\sum_{i=1}^{n} \psi(X_i;\theta_n)\right| \stackrel{P_{\theta_0}}{\to} 0 \tag{2.1}$$

for all sequences $\{\theta_n\}$ with $\{\sqrt{n}(\theta_n - \theta_0)\}$ bounded. We call $\psi$ the *influence function* of $\hat{\kappa}_n$.

Suppose we have an estimator $\hat{\theta}_n = t_n(X_1, \ldots, X_n)$, $t_n : \mathcal{X}^n \to \mathbb{R}^k$, $\mathcal{A}^n$-Borel measurable, that is a locally asymptotically linear estimator of $\theta$ at $\theta_0$ with influence function $\psi_\theta \in \Psi_k$, that is,

$$\sqrt{n}\left|\hat{\theta}_n - \theta_n - \frac{1}{n}\sum_{i=1}^{n}\psi_\theta(X_i;\theta_n)\right| \stackrel{P_{\theta_0}}{\to} 0 \tag{2.2}$$

holds for all sequences $\{\theta_n\}$ with $\{\sqrt{n}(\theta_n - \theta_0)\}$ bounded. Suppose furthermore that there is a process $\hat{\kappa}_{\theta,n} = k_n(X_1, \ldots, X_n; \theta)$ that is locally asymptotically linear in $\psi_\kappa \in \Psi_m$ around $\kappa(\theta)$ such that

$$(2.3) \qquad \sqrt{n}\left|\hat{\kappa}_{\theta_n,n} - \kappa(\theta_n) - n^{-1}\sum_{i=1}^n \psi_\kappa(X_i; \theta_n)\right| \stackrel{P_{\theta_0}}{\to} 0$$

holds for all sequences $\{\theta_n\}$ with $\{\sqrt{n}(\theta_n - \theta_0)\}$ bounded. Note that we have extended here the concept of local asymptotic linearity from estimators of $\kappa(\theta)$, as in Definition 2.1, to statistics indexed by $\theta$. This is quite reasonable since the gist of the concept is that the relevant statistic behaves as an average locally asymptotically.

We want to describe the local asymptotic behavior of the substitution estimator

$$(2.4) \qquad \hat{\kappa}_{n,1} = \hat{\kappa}_{\hat{\theta}_n,n},$$

which replaces the unknown $\theta$ by its estimator $\hat{\theta}_n$ in $\hat{\kappa}_{\theta,n}$.

Heuristically, by (2.3), $\sqrt{n}(\hat{\kappa}_{\hat{\theta}_n,n} - \kappa(\theta_0))$ behaves like

$$(2.5) \qquad \sqrt{n}\left(\kappa(\hat{\theta}_n) - \kappa(\theta_0) + n^{-1}\sum_{i=1}^n \psi_\kappa(X_i; \hat{\theta}_n)\right).$$

Now it is natural to assume the existence of a matrix-valued function $c: \Theta \to \mathbb{R}^{m \times k}$ that is continuous at $\theta_0$ and that is such that for every sequence $\{\theta_n\}$ with $\{\sqrt{n}(\theta_n - \theta)\}$ bounded,

$$(2.6) \quad \left|\frac{1}{\sqrt{n}}\sum_{i=1}^n \psi_\kappa(X_i; \theta_n) - \frac{1}{\sqrt{n}}\sum_{i=1}^n \psi_\kappa(X_i; \theta_0) - c(\theta_0)\sqrt{n}(\theta_n - \theta_0)\right| \stackrel{P_{\theta_0}}{\to} 0$$

holds. Since $\kappa$ is differentiable, $\sqrt{n}(\hat{\kappa}_{\hat{\theta}_n,n} - \kappa(\theta_0))$ would then behave like

$$\sqrt{n}\left(\kappa'(\theta_0)(\hat{\theta}_n - \theta_0) + n^{-1}\sum_{i=1}^n \psi_\kappa(X_i; \theta_0) + c(\theta_0)(\hat{\theta}_n - \theta_0)\right)$$

and hence, by (2.2), like

$$\frac{1}{\sqrt{n}}\sum_{i=1}^n (\psi_\kappa(X_i; \theta_0) + (\kappa'(\theta_0) + c(\theta_0))\psi_\theta(X_i; \theta_0)).$$

The estimator $\hat{\kappa}_{\hat{\theta}_n,n}$ thus inherits its asymptotic linearity from the submodel estimator $\hat{\kappa}_{\theta,n}$ and the estimator $\hat{\theta}_n$. To study this asymptotic linearity more carefully, we first describe a sample splitting procedure, for which we can prove statements under minimal conditions. Fix a sequence of integers $\{\lambda_n\}_{n=1}^\infty$, such that

$$(2.7) \qquad \frac{\lambda_n}{n} \to \frac{1}{2}.$$

We split the sample $(X_1,\ldots,X_n)$ into two parts, $(X_1,\ldots,X_{\lambda_n})$ and $(X_{\lambda_n+1},\ldots,X_n)$. Define

$$(2.8) \quad \begin{aligned} \tilde{\theta}_{n1} &= t_{\lambda_n}(X_1,\ldots,X_{\lambda_n}), & \tilde{\theta}_{n2} &= t_{n-\lambda_n}(X_{\lambda_n+1},\ldots,X_n), \\ \tilde{\kappa}^{(1)}_{\theta,\lambda_n} &= k_{\lambda_n}(X_1,\ldots,X_{\lambda_n};\theta), & \tilde{\kappa}^{(2)}_{\theta,n-\lambda_n} &= k_{n-\lambda_n}(X_{\lambda_n+1},\ldots,X_n;\theta) \end{aligned}$$

and

$$(2.9) \quad \hat{\kappa}_{n,2} = \frac{\lambda_n}{n}\tilde{\kappa}^{(1)}_{\tilde{\theta}_{n2},\lambda_n} + \frac{n-\lambda_n}{n}\tilde{\kappa}^{(2)}_{\tilde{\theta}_{n1},n-\lambda_n}.$$

The following theorem describes the influence function of this split-sample substitution estimator $\hat{\kappa}_{n,2}$.

THEOREM 2.1. *Fix $\theta_0 \in \Theta$. Suppose that $\kappa:\Theta \to \mathbb{R}^m$ is continuously differentiable in a neighborhood of $\theta_0$ with derivative matrix $\kappa'$ and suppose that conditions* (2.2), (2.3) *and* (2.6) *hold for some $c:\Theta \to \mathbb{R}^{m\times k}$ that is continuous at $\theta_0$. Then $\hat{\kappa}_{n,2}$ defined by* (2.7)–(2.9) *is locally asymptotically linear for $\kappa$ at $\theta_0$ with influence function $\tilde{\psi}$ given by*

$$(2.10) \quad \tilde{\psi}(x;\theta) = \psi_\kappa(x;\theta) + (\kappa'(\theta) + c(\theta))\psi_\theta(x;\theta),$$

*that is, for every sequence $\{\theta_n\}$ with $\{\sqrt{n}(\theta_n - \theta_0)\}$ bounded,*

$$(2.11) \quad \sqrt{n}\left|\hat{\kappa}_{n,2} - \kappa(\theta_n) - \frac{1}{n}\sum_{i=1}^n \tilde{\psi}(X_i;\theta_n)\right| \overset{P_{\theta_0}}{\to} 0.$$

PROOF. Fix $\theta_0 \in \Theta$ and the sequence $\{\theta_n\}$. Take another sequence $\{\tilde{\theta}_n\}$ such that $\{\sqrt{n}(\tilde{\theta}_n - \theta_0)\}$ stays bounded. Combining (2.3), with $\theta_n$ replaced by $\tilde{\theta}_n$, and (2.6), both with $\theta_n$ and with $\theta_n$ replaced by $\tilde{\theta}_n$, we obtain, using $(X_1,\ldots,X_{\lambda_n})$,

$$\sqrt{n}\left|\tilde{\kappa}^{(1)}_{\tilde{\theta}_n,\lambda_n} - \kappa(\tilde{\theta}_n) - \frac{1}{\lambda_n}\sum_{i=1}^{\lambda_n}\psi_\kappa(X_i;\theta_n) - c(\theta_0)(\tilde{\theta}_n - \theta_n)\right| \overset{P_{\theta_0}}{\to} 0,$$

which by continuous differentiability of $\kappa(\cdot)$ and continuity of $c(\cdot)$ at $\theta_0$ yields

$$(2.12) \quad \sqrt{n}\left|\tilde{\kappa}^{(1)}_{\tilde{\theta}_n,\lambda_n} - \kappa(\theta_n) - \frac{1}{\lambda_n}\sum_{i=1}^{\lambda_n}\psi_\kappa(X_i;\theta_n) - (\kappa'(\theta_n) + c(\theta_n))(\tilde{\theta}_n - \theta_n)\right| \overset{P_{\theta_0}}{\to} 0.$$

By the asymptotic linearity of $\hat{\theta}_n$, we have

$$(2.13) \quad \sqrt{n}\left|\tilde{\theta}_{n2} - \theta_n - \frac{1}{n-\lambda_n}\sum_{i=\lambda_n+1}^n \psi_\theta(X_i;\theta_n)\right| \overset{P_{\theta_0}}{\to} 0.$$

Hence, by the independence of $(X_1, \ldots, X_{\lambda_n})$ and $(X_{\lambda_n+1}, \ldots, X_n)$, (2.12) and (2.13) together yield

$$\sqrt{n}\left|\tilde{\kappa}^{(1)}_{\tilde{\theta}_{n2},\lambda_n} - \kappa(\theta_n) - \frac{1}{\lambda_n}\sum_{i=1}^{\lambda_n}\psi_\kappa(X_i;\theta_n)\right.$$
$$\left. - \frac{1}{n-\lambda_n}\sum_{i=\lambda_n+1}^{n}(\kappa'(\theta_n) + c(\theta_n))\psi_\theta(X_i;\theta_n)\right| \overset{P_{\theta_0}}{\to} 0.$$

Similarly we obtain

$$\sqrt{n}\left|\tilde{\kappa}^{(2)}_{\tilde{\theta}_{n1},n-\lambda_n} - \kappa(\theta_n) - \frac{1}{n-\lambda_n}\sum_{i=\lambda_n+1}^{n}\psi_\kappa(X_i;\theta_n)\right.$$
$$\left. - \frac{1}{\lambda_n}\sum_{i=1}^{\lambda_n}(\kappa'(\theta_n) + c(\theta_n))\psi_\theta(X_i;\theta_n)\right| \overset{P_{\theta_0}}{\to} 0.$$

These last two statements yield

$$\sqrt{n}\left|\hat{\kappa}_{n,2} - \kappa(\theta_n) - \frac{1}{n}\sum_{i=1}^{n}\psi_\kappa(X_i;\theta_n)\right.$$
(2.14)
$$- \frac{1}{n}\sum_{i=1}^{n}\left\{\frac{n-\lambda_n}{\lambda_n}\mathbf{1}_{[i\leq\lambda_n]} + \frac{\lambda_n}{n-\lambda_n}\mathbf{1}_{[i>\lambda_n]}\right\}$$
$$\left. \times (\kappa'(\theta_n) + c(\theta_n))\psi_\theta(X_i;\theta_n)\right| \overset{P_{\theta_0}}{\to} 0.$$

In view of (2.7) this shows that $\hat{\kappa}_{n,2}$ is a locally asymptotically linear estimator of $\kappa$, with influence function given by (2.10). This proves the theorem. □

Note that the expression within braces in (2.14) reveals why (2.7) is crucial to our sample splitting scheme.

To establish local asymptotic linearity of the direct substitution estimator $\hat{\kappa}_{n,1}$ without sample splitting [cf. (2.4)], we need locally asymptotically uniform continuity in $\theta$ at $\theta_0$ of the estimators $\hat{\kappa}_{\theta,n}$ as follows:

For every $\delta > 0$, $\varepsilon > 0$ and $c > 0$, there exist $\zeta > 0$ and $n_0 \in \mathbb{N}$ such that for all $n \geq n_0$

(2.15) $$P_{\theta_0}\left(\sup_{\sqrt{n}|\theta-\theta_0|\leq c, \sqrt{n}|\theta-\tilde{\theta}|\leq \zeta}\sqrt{n}|\hat{\kappa}_{\theta,n} - \hat{\kappa}_{\tilde{\theta},n}| \geq \varepsilon\right) \leq \delta.$$

THEOREM 2.2. *Fix $\theta_0 \in \Theta$. If (2.15) holds in the situation of Theorem 2.1, then the substitution estimator $\hat{\kappa}_{n,1} = \hat{\kappa}_{\hat{\theta}_n,n}$ is a locally asymptotically linear estimator of $\kappa$ with influence function $\tilde{\psi}$ given by (2.10).*

PROOF. Fix $\theta_0 \in \Theta$, $\delta > 0$, $\varepsilon > 0$. Choose $c$ and $n_0$ such that for $n \geq n_0$

$$(2.16) \qquad P_{\theta_0}(\sqrt{n}|\hat{\theta}_n - \theta_0| > c) \leq \delta.$$

Now, choose $\zeta$ sufficiently small such that (2.15) holds too (increase $n_0$ if necessary), and such that the matrix norm of $\kappa'(\theta_0) + c(\theta_0)$ satisfies

$$(2.17) \qquad \|\kappa'(\theta_0) + c(\theta_0)\|\zeta < \varepsilon/2.$$

Let $\hat{\theta}_n(\zeta)$ be the efficient estimator $\hat{\theta}_n$ discretized via a grid $\mathcal{G}_\zeta$ of meshwidth $2(kn)^{-1/2}\zeta$, such that

$$(2.18) \qquad \sqrt{n}|\hat{\theta}_n(\zeta) - \hat{\theta}_n| \leq \zeta \qquad \text{a.s.}$$

It follows by (2.18) and (2.15) that for $n \geq n_0$ the inequality

$$(2.19) \begin{aligned}
& P_{\theta_0}\Bigg(\sqrt{n}\bigg|\hat{\kappa}_{\hat{\theta}_n,n} - \kappa(\theta_n) \\
&\qquad - \frac{1}{n}\sum_{i=1}^n [\psi_\kappa(X_i;\theta_n) + (\kappa'(\theta_n) + c(\theta_n))\psi_\theta(X_i;\theta_n)]\bigg| \geq 4\varepsilon\Bigg) \\
&= P_{\theta_0}\Bigg(\sqrt{n}\bigg|(\hat{\kappa}_{\hat{\theta}_n,n} - \hat{\kappa}_{\hat{\theta}_n(\zeta),n}) \\
&\qquad + \bigg\{\hat{\kappa}_{\hat{\theta}_n(\zeta),n} - \kappa(\theta_n) - \frac{1}{n}\sum_{i=1}^n \psi_\kappa(X_i;\theta_n) \\
&\qquad\qquad - (\kappa'(\theta_n) + c(\theta_n))(\hat{\theta}_n(\zeta) - \theta_n)\bigg\} \\
&\qquad + (\kappa'(\theta_n) + c(\theta_n))(\hat{\theta}_n(\zeta) - \hat{\theta}_n) \\
&\qquad + (\kappa'(\theta_n) + c(\theta_n))\bigg\{\hat{\theta}_n - \theta_n - \frac{1}{n}\sum_{i=1}^n \psi_\theta(X_i;\theta_n)\bigg\}\bigg| \geq 4\varepsilon\Bigg) \\
&\leq P_{\theta_0}(\sqrt{n}|\hat{\theta}_n - \theta_0| > c) + \delta \\
&\qquad + \sum_{\tilde{\theta} \in \mathcal{G}_\zeta, \sqrt{n}|\tilde{\theta} - \theta_0| \leq c+\zeta} P_{\theta_0}\Bigg(\sqrt{n}\bigg|\hat{\kappa}_{\tilde{\theta},n} - \kappa(\theta_n) - \frac{1}{n}\sum_{i=1}^n \psi_\kappa(X_i;\theta_n) \\
&\qquad\qquad - (\kappa'(\theta_n) + c(\theta_n))(\tilde{\theta} - \theta_n)\bigg| \geq \varepsilon\Bigg)
\end{aligned}$$

$$+ P_{\theta_0}(\|\kappa'(\theta_n) + c(\theta_n)\|\zeta \geq \varepsilon)$$
$$+ P_{\theta_0}\left(\|\kappa'(\theta_n) + c(\theta_n)\|\sqrt{n}\left|\hat{\theta}_n - \theta_n - \frac{1}{n}\sum_{i=1}^n \psi_\theta(X_i;\theta_n)\right| \geq \varepsilon\right)$$

holds. In view of (2.16), in view of the boundedness of the number of terms in the sum with all terms converging to zero by (2.12), in view of (2.17) and the continuity of $\kappa' + c$, and in view of the linearity of $\hat{\theta}_n$ [see (2.2)], the lim sup as $n \to \infty$ of the right-hand side of (2.19) equals at most $2\delta$. Since $\delta$ may be chosen arbitrarily small, this proves the asymptotic linearity. □

REMARK 2.1. In some cases, it may happen that $\kappa'(\cdot) + c(\cdot) = 0$. Then it is easily seen that the influence function of the substitution estimators $\hat{\kappa}_{n,1}$ and $\hat{\kappa}_{n,2}$ is given by $\tilde{\psi}(\cdot;\cdot) = \psi_\kappa(\cdot;\cdot)$, even if $\hat{\theta}_n$ is not locally asymptotically linear but is just $\sqrt{n}$-consistent.

REMARK 2.2. If $\psi_\kappa(\cdot;\theta)$ is differentiable in $\theta$ with derivative $\dot{\psi}_\kappa(\cdot;\theta)$, then Taylor expansion and the law of large numbers suggest $c(\theta) = E_\theta \dot{\psi}_\kappa(X_1;\theta)$. Furthermore, differentiation of $E_\theta \psi_\kappa(X_1;\theta) = 0$ with respect to $\theta$ hints at

$$(2.20) \qquad c(\theta) = E_\theta \dot{\psi}_\kappa(X_1;\theta) = -E_\theta \psi_\kappa(X_1;\theta)\dot{l}^\top(X_1;\theta),$$

with $\dot{l}(x;\theta)$ the score function for $\theta$, namely $\partial \log p(x;\theta)/\partial \theta$.

REMARK 2.3. Theorem 2.2 is related to a result known as the delta method; see Section 2.5 of Lehmann (1999). Given the function $\kappa(\cdot)$, choose $\hat{\kappa}_{\theta,n} = \kappa(\theta)$. Then the convergence (2.3) holds trivially with $\psi_\kappa(\cdot;\cdot) = 0$. Furthermore, (2.6) is valid with $c(\cdot) = 0$ and (2.15) holds if $\kappa$ is continuously differentiable. Now Theorem 2.2 states that the local asymptotic linearity of $\hat{\theta}_n$ in (2.2) implies the local asymptotic linearity of $\kappa(\hat{\theta}_n)$, that is,

$$\sqrt{n}\left|\kappa(\hat{\theta}_n) - \kappa(\theta_n) - \frac{1}{n}\sum_{i=1}^n \kappa'(\theta_n)\psi_\theta(X_i;\theta_n)\right| \stackrel{P_{\theta_0}}{\to} 0,$$

and hence by the central limit theorem the asymptotic normality of $\sqrt{n}(\kappa(\hat{\theta}_n) - \kappa(\theta_0))$ under $P_{\theta_0}$. Note that the delta method states that asymptotic normality of $\sqrt{n}(\hat{\theta}_n - \theta_0)$ implies asymptotic normality of $\sqrt{n}(\kappa(\hat{\theta}_n) - \kappa(\theta_0))$.

REMARK 2.4. Under different sets of regularity conditions, the heredity of asymptotic normality of substitution statistics has been proved by Randles (1982) and Pierce (1982). Since asymptotic linearity implies asymptotic normality, both our conditions and our conclusions in proving heredity of asymptotic linearity are stronger than needed for heredity of asymptotic normality. However, the approach via differentiability in Section 3 of Randles (1982) comes pretty close to the assumption of asymptotic linearity. Moreover, our ultimate goal is the study of efficient estimators, which are bound to be asymptotically linear.

Let us now discuss sufficient conditions for (2.6). The following standard result will be quite helpful and may be verified by studying first and second moments.

LEMMA 2.1. *Fix $\theta_0 \in \Theta$. If*

$$(2.21) \qquad E_{\theta_0}(|\psi_\kappa(X_1; \theta_0 + \varepsilon) - \psi_\kappa(X_1; \theta_0)|^2) \to 0 \qquad as\ \varepsilon \to 0$$

*holds and the map $\varepsilon \mapsto E_{\theta_0} \psi_\kappa(X_1; \theta_0 + \varepsilon)$ is differentiable at 0 with derivative matrix $c(\theta_0)$, then (2.6) holds for all sequences $\{\theta_n\}$ with $\{\sqrt{n}(\theta_n - \theta_0)\}$ bounded.*

Sometimes (2.6) may be verified by a direct application of the following "law-of-large-numbers"-type of result.

LEMMA 2.2. *Fix $\theta_0 \in \Theta$. For $P_{\theta_0}$-almost all $x$, let $\psi_\kappa(x; \theta)$ be continuously differentiable in $\theta$ with derivative $\dot{\psi}_\kappa(x; \theta)$. If*

$$(2.22) \qquad E_{\theta_0} |\dot{\psi}_\kappa(X_1; \theta)| \to E_{\theta_0} |\dot{\psi}_\kappa(X_1; \theta_0)| \qquad as\ \theta \to \theta_0,$$

*then for $\{\sqrt{n}(\theta_n - \theta_0)\}$ bounded*

$$(2.23)\quad \frac{1}{\sqrt{n}} \sum_{i=1}^n \{\psi_\kappa(X_i; \theta_n) - \psi_\kappa(X_i; \theta_0)\} - \sqrt{n}(\theta_n - \theta_0)^\top E_{\theta_0} \dot{\psi}_\kappa(X_1; \theta_0) \overset{P_{\theta_0}}{\to} 0$$

*holds.*

PROOF. Write the left-hand side of (2.23) as

$$\sqrt{n}(\theta_n - \theta_0)^\top \left[ \frac{1}{n} \sum_{i=1}^n \{\dot{\psi}_\kappa(X_i; \theta_0) - E_{\theta_0} \dot{\psi}_\kappa(X_i; \theta_0)\} \right.$$
$$\left. + \int_0^1 \frac{1}{n} \sum_{i=1}^n \{\dot{\psi}_\kappa(X_i; \theta_0 + \zeta(\theta_n - \theta_0)) - \dot{\psi}_\kappa(X_i; \theta_0)\} \, d\zeta \right],$$

note that the first absolute moment of the last term may be bounded by

$$\sqrt{n}|\theta_n - \theta_0| \int_0^1 E_{\theta_0} |\dot{\psi}_\kappa(X_1; \theta_0 + \zeta(\theta_n - \theta_0)) - \dot{\psi}_\kappa(X_1; \theta_0)| \, d\zeta$$

and apply, for example, Theorem A.7.2 of Bickel, Klaassen, Ritov and Wellner (1993). □

Condition (2.6) may be also derived via regularity of $\mathcal{P}$ and local asymptotic normality (LAN) by an argument similar to the one leading to (2.1.15) of Proposition 2.1.2, pages 16 and 17, of Bickel, Klaassen, Ritov and Wellner (1993).

DEFINITION 2.2. A parametric model $\mathcal{P} = \{P_\theta : \theta \in \Theta\}$, $\Theta \subset \mathbb{R}^k$ open, is a $k$-dimensional *regular parametric model* if there exists a $\sigma$-finite dominating measure $\mu$ such that, with $p(\theta) = dP_\theta/d\mu$, $s(\theta) = p_\theta^{1/2}$:

(i) for all $\theta \in \Theta$ there exists a $k$-vector $\dot{l}(\theta)$ of score functions in $L_2(P_\theta)$ such that

$$(2.24) \qquad s(\tilde{\theta}) = s(\theta) + \tfrac{1}{2}(\tilde{\theta} - \theta)^\top \dot{l}(\theta) s(\theta) + o(|\tilde{\theta} - \theta|)$$

in $L_2(\mu)$ as $|\tilde{\theta} - \theta| \to 0$;

(ii) for every $\theta \in \Theta$ the $k \times k$ Fisher information matrix $\int \dot{l}(\theta) \dot{l}^\top(\theta) p(\theta) \, d\mu$ is nonsingular;

(iii) the map $\theta \mapsto \dot{l}(\theta) s(\theta)$ is continuous from $\Theta$ to $L_2^k(\mu)$.

A priori, it would have been more general if condition (i) of Definition 2.2 had prescribed Fréchet-differentiability of $s(\theta)$ with derivative $\dot{s}(\theta)$ in $L_2^k(\mu)$. However, it can be shown that all components of $\dot{s}(\theta)$ would vanish then almost everywhere where $s(\theta)$ vanishes. Consequently, $\dot{s}(\theta)$ may be written as $\dot{l}(\theta) s(\theta)/2$ in $L_2^k(\mu)$; see Proposition A.5.3.F of Bickel, Klaassen, Ritov and Wellner (1993).

This approach to prove (2.6) through regularity and local asymptotic normality has been implemented in a preprint of the present paper [Klaassen and Putter (2000)]. However, a much nicer argument has been noted by Schick (2001).

LEMMA 2.3. *Suppose that the model $\mathcal{P}$ is regular and fix $\theta_0 \in \Theta$. If $\psi_\kappa \in \Psi_m$ satisfies the continuity condition*

$$(2.25) \qquad \|\psi_\kappa(\cdot; \tilde{\theta}) s(\tilde{\theta}) - \psi_\kappa(\cdot; \theta) s(\theta)\|_\mu \to 0,$$

*as $\tilde{\theta} \to \theta$, then (2.6) is valid with $c(\theta)$ given by (2.20).*

PROOF. Since the regularity of $\mathcal{P}$ implies Hellinger differentiability at $\theta_0$, Theorem 2.3 of Schick (2001) may be applied and yields (2.6). The continuity of $c(\cdot)$ is implied by (2.25) and the regularity of $\mathcal{P}$. □

REMARK 2.5. At the end of his Section 1 on page 17, Schick (2001) refers to (3.5) of the preprint Klaassen and Putter (2000). This is just (2.6) of the present version of this paper.

**3. Examples for asymptotic linearity of substitution estimators.** Although the results of Section 2 are stated within a parametric model, most of the applications we have in mind (in particular efficiency as discussed in Section 5) are in the context of semiparametric models where the interest is in a functional of the infinite-dimensional parameter only. In the analysis of these applications it suffices to study parametric submodels where the infinite-dimensional parameter is fixed. Hence the results of Section 2 are also applicable in this context. In order to illustrate the heredity of asymptotic linearity of substitution estimators in the framework of semiparametric models, however, we need to introduce some notation and conventions specific to semiparametric models.

Let $\mathcal{P} = \{P_{\theta,G} : \theta \in \Theta, G \in \mathcal{G}\}$, $\Theta \subset \mathbb{R}^k$ open, $\mathcal{G} \subset \mathcal{H}$, be our semiparametric model (1.1). The model $\mathcal{P}$ might be parametric in the sense that $\mathcal{G}$ is Euclidean. We may represent the elements of $\mathcal{P}$ by the square roots $s(\theta, G) = p^{1/2}(\theta, G)$ of their densities $p(\theta, G)$ with respect to a $\sigma$-finite dominating measure $\mu$ if such a dominating measure exists on the sample space $(\mathcal{X}, \mathcal{A})$.

By keeping $G$ fixed and by varying $\theta$ over $\Theta$ we get a parametric submodel of $\mathcal{P}$, denoted by $\mathcal{P}_1 = \mathcal{P}_1(G)$. Often $\mathcal{P}_1(G)$ will be a regular parametric model in the sense of Definition 2.2. The $k$-vector of score functions of $\mathcal{P}_1(G)$ will be denoted by $\dot{l}_1$ then and in particular we will have

$$(3.1) \qquad s(\tilde{\theta}, G) = s(\theta, G) + \tfrac{1}{2}(\tilde{\theta} - \theta)^\top \dot{l}_1(\theta, G) s(\theta, G) + \mathcal{o}(|\tilde{\theta} - \theta|).$$

Let $X_1, \ldots, X_n$ be an i.i.d. sample from $P_{\theta,G} \in \mathcal{P}$ and let $\kappa : \mathcal{P} \to \mathbb{R}^m$ be an unknown Euclidean parameter of the model $\mathcal{P}$ with

$$(3.2) \qquad \kappa(P_{\tilde{\theta},G}) = \kappa(P_{\theta,G}) = \tilde{\kappa}(G), \qquad \theta, \tilde{\theta} \in \Theta, G \in \mathcal{G},$$

for some $\tilde{\kappa} : \mathcal{G} \to \mathbb{R}^m$. Since interest is mainly in estimating a Banach parameter $\nu = \tilde{\nu}(G) \in \mathcal{B}$ as in (1.2), a typical choice of $\kappa$ with $m = 1$ would be $\kappa(P_{\theta,G}) = b^* \tilde{\nu}(G)$ for some $b^* \in \mathcal{B}^*$, the dual of $\mathcal{B}$; note that such a parameter $\kappa$ is independent of $\theta$ in the sense of (3.2). Let $\hat{\kappa}_n = k_n(X_1, \ldots, X_n)$ be an estimator of $\kappa$ with $k_n : \mathcal{X}^n \to \mathbb{R}^m$ an $\mathcal{A}^n$-Borel measurable function. As in Definition 2.1, the estimator $\hat{\kappa}_n$ of $\kappa(P_{\theta,G}) = \tilde{\kappa}(G)$ is called *locally asymptotically linear* at $P_{\theta_0,G}$ if there exists a measurable function $\psi(\cdot; \cdot, G) \in \Psi_m$ such that

$$(3.3) \qquad \sqrt{n} \left| \hat{\kappa}_n - \tilde{\kappa}(G) - \frac{1}{n} \sum_{i=1}^n \psi(X_i; \theta_n, G) \right| \stackrel{P_{\theta_0,G}}{\to} 0$$

holds for all sequences $\{\theta_n\}$ with $\{\sqrt{n}(\theta_n - \theta_0)\}$ bounded. The function $\psi(\cdot; \theta, G)$ is called the *influence function* of $\hat{\kappa}_n$ at $P_{\theta,G}$ and $\psi(\cdot; \cdot, G)$ is called influence function as well.

The results of Section 2 are illustrated in the following examples.

EXAMPLE 3.1 (Sample variance). Let $X_1,\ldots,X_n$ be i.i.d. with distribution function $G(\cdot - \theta)$ on $\mathbb{R}$. Here $G$ is an unknown distribution function with mean zero and finite fourth moment. Given $\theta$, a good estimator of the variance $\tilde{\kappa}(G) = \int x^2\, dG(x)$ of $G$ is $\hat{\kappa}_{\theta,n} = n^{-1}\sum_{i=1}^n (X_i - \theta)^2$, which is linear with influence function

$$\psi_\kappa(x;\theta,G) = (x-\theta)^2 - \tilde{\kappa}(G). \tag{3.4}$$

Since $\theta$ can be estimated by the sample mean $\hat{\theta}_n = \bar{X}_n$, which is linear and hence asymptotically linear, Theorem 2.2 yields the sample variance

$$\hat{\kappa}_{\hat{\theta}_n,n} = S_n^2 = \frac{1}{n}\sum_{i=1}^n (X_i - \bar{X}_n)^2 \tag{3.5}$$

as a locally asymptotically linear estimator of $\tilde{\kappa}(G)$ in case $\theta$ is unknown; note that (2.15) holds in view of the law of large numbers. The sample variance is adaptive in the sense that it has the same influence function $\psi_\kappa$ as in (3.4) because (2.6) holds with $c(\theta) = 0$, as may be verified easily.

Of course this estimator is the prototype of a substitution estimator, used routinely to the extent that typically it is not recognized as a substitution estimator.

EXAMPLE 3.2 (Parametrized linkage models). Observe realizations of $X_i$, $i=1,\ldots,n$, that are i.i.d. copies of $X$. In many statistical models the random variable $X$ is linked to an error random variable $\epsilon$ with distribution function $G$. This linkage is parametrized by $\theta \in \Theta \subset \mathbb{R}^k$ and may be described by a measurable map $t_\theta : \mathcal{X} \to \mathbb{R}$ with

$$t_\theta(X) = \epsilon.$$

The prime example is the linear regression model with

$$t_\theta(x) = y - \theta^\top z, \qquad x = (y, z^\top)^\top, y \in \mathbb{R}, z \in \mathbb{R}^k, E\epsilon = 0,$$

yielding the error random variable $\epsilon$ and

$$t_\theta(x) = (y - \nu^\top z)/\sigma, \qquad \theta = (\nu^\top, \sigma)^\top, x = (y, z^\top)^\top,\ y \in \mathbb{R},\ z \in \mathbb{R}^{k-1},$$

with $E\epsilon = 0$, $E\epsilon^2 = Et_\theta^2(X) = 1$, generating the standardized random variable $\epsilon$. Another example is the accelerated failure time model with

$$t_\theta(x) = e^{-\theta^\top z} y, \qquad x = (y, z^\top)^\top, y \in [0, \infty),\ z \in \mathbb{R}^k,$$

yielding the standardized life time random variable $\epsilon$. Recall that the distribution of $X$ is denoted by $P_{\theta,G}$. The Euclidean parameter is $\theta$ and the error distribution function or the standardized life time distribution function $\nu(P_{\theta,G}) = G$ could be the Banach parameter of interest. Given $\theta$, an obvious estimator of $G$ would be the empirical distribution function of $t_\theta(X_i)$, $i = 1,\ldots, n$.

We will study estimation of the one-dimensional parameter

$$\kappa(P_{\theta,G}) = \tilde{\kappa}(G) = \int h\,dG,$$

where $h$ is some known function with $\int h^2\,dG < \infty$. Taking the empirical distribution function of $t_\theta(X_i)$ as an estimator of $G$ when $\theta$ is known, we obtain

$$\hat{\kappa}_{\theta,n} = \frac{1}{n}\sum_{i=1}^n h(t_\theta(X_i))$$

as an estimator of $\kappa(P_{\theta,G})$. This estimator is linear and hence locally asymptotically linear in the sense of (2.3) in the influence function

$$\psi_\kappa(x;\theta,G) = h(t_\theta(x)) - \kappa(P_{\theta,G}), \qquad x \in \mathcal{X}.$$

If $\hat{\theta}_n$ is a locally asymptotically linear estimator of $\theta$ with influence function $\psi_\theta(\cdot;\theta,G) \in \Psi_k$ as in (2.2), an application of Lemma 2.1 yields the validity of Theorem 2.1 provided

$$E_{\theta_0,G}(|h(t_{\theta_0+\varepsilon}(X)) - h(t_{\theta_0}(X))|^2) \to 0 \qquad \text{as } \varepsilon \to 0$$

holds and $E_{\theta_0,G}h(t_{\theta_0+\varepsilon}(X))$ is differentiable in $\varepsilon$ at 0 with a derivative matrix $c(\theta_0)$ that is continuous in $\theta_0$. Noting that $\kappa'(\theta)$ from (2.10) vanishes here, we see that the local asymptotic linearity of the split-sample substitution estimator $\hat{\kappa}_{n,2}$ from (2.9) holds with influence function

$$\tilde{\psi}(x;\theta,G) = h(t_\theta(x)) - \kappa(P_{\theta,G}) + c(\theta)\psi_\theta(x;\theta,G).$$

Note that the sample variance is a special case with $h(x) = x^2$, $t_\theta(x) = x - \theta$ and $c(\theta) = 0$.

In Section 6 we shall consider the most important special case of this example, the linear regression model, in more detail in the context of efficient estimation of (functionals of the) error distribution. Here we consider the linear regression model with standardized errors. Substitution of, for example, the least squares estimators would lead to the empirical distribution of the standardized residuals as a natural estimator of $G$. See, for example, Koul (1992, 2002) or Loynes (1980) for early studies of the empirical distribution function of regression residuals. With $\kappa(P) = \tilde{\kappa}(G) = \int h\,dG = Eh(\epsilon)$ for appropriate functions $h: \mathbb{R} \to \mathbb{R}$, Theorems 2.1 and 2.2 hold with

$$\psi_\kappa(x;\theta,G) = h\left(\frac{y - \alpha - \beta z}{\sigma}\right) - \int h\,dG$$

and

$$c(\theta) = -\frac{1}{\sigma}(Eh'(\epsilon), Eh'(\epsilon)EZ, E\varepsilon h'(\epsilon))^\top.$$

For $h(\epsilon) = \epsilon^3$ we obtain an estimator of the skewness of the error distribution, whose asymptotic normality has been studied by Pierce (1982) under normality (see also Remark 2.4).

EXAMPLE 3.3 (Distribution function in two-sample location model). $X_1, \ldots, X_n$ are i.i.d. copies of $X = (Y, Z)$, where $Y$ and $Z - \theta$ are i.i.d. with density $g$, $\int y^2 g(y) \, dy < \infty$. The Banach parameter of interest is the distribution function $G(\cdot) = \int_{-\infty}^{\cdot} g(y) \, dy$. Given the shift parameter $\theta$, it can be estimated by the linear estimator

$$(3.6) \qquad \hat{G}_{\theta,n}(y) = \frac{1}{2n} \sum_{i=1}^{n} (\mathbf{1}_{[Y_i \leq y]} + \mathbf{1}_{[Z_i - \theta \leq y]}), \qquad y \in \mathbb{R}.$$

Since $\hat{\theta}_n = \bar{Z}_n - \bar{Y}_n$ is a linear estimator of $\theta$, the substitution estimator $\hat{G}_{\hat{\theta}_n,n}(\cdot)$ is locally asymptotically linear by Theorem 2.2. A structure similar to the one in Example 3.2 may be described by the map $t_\theta : \mathcal{X} \to \mathbb{R}^2$ with $t_\theta(x) = (y, z - \theta)^\top$, $y, z \in \mathbb{R}$.

In fact, $\theta$ can be estimated adaptively in Example 3.3, that is, efficiently within this semiparametric model; see van Eeden (1970) for an early construction valid for the class of strongly unimodal densities $g$ and Beran (1974) and Stone (1975) for the most general situation. If we apply such an asymptotically efficient estimator $\hat{\theta}_n$, then the resulting estimator $\hat{G}_{\hat{\theta}_n,n}(\cdot)$ is asymptotically efficient too, since (3.6) is efficient given $\theta$. This hereditary property of asymptotic efficiency for substitution estimators follows from the heredity for linearity, which will be shown in Section 5 and is the main result of the present paper. As preparation we study efficient influence functions in the next section.

**4. Efficient influence functions.** Let $\mathcal{H}$ be a Hilbert space. A one-dimensional subset $\{h_\eta \in \mathcal{H} : -1 < \eta < 1\}$ of $\mathcal{H}$ is called a *path* if the map $\eta \mapsto h_\eta$ is continuously Fréchet differentiable with nonvanishing derivative, implying, for example, the existence of an $\dot{h} \in \mathcal{H}$, $\dot{h} \neq 0$, with

$$(4.1) \qquad h_\eta = h_0 + \eta \dot{h} + \mathcal{o}(\eta) \qquad \text{in } \mathcal{H}$$

as $\eta \to 0$.

Let $\mathcal{P}$ be a statistical model, that is, a collection of probability distributions, and fix $P \in \mathcal{P}$. A subset $\mathcal{P}_P = \{P_\eta : -1 < \eta < 1\}$ of $\mathcal{P}$ is called a path through $P$ if $P_0$ equals $P$ and $\mathcal{P}_P$ is a regular one-dimensional parametric submodel in the sense of Definition 2.2. This implies the existence of a so-called tangent $t \in L_2(P)$, $t \neq 0$, such that with $s_\eta = \sqrt{dP_\eta/d\mu}$ for some dominating $\sigma$-finite measure $\mu$, and with $s = s_0$,

$$(4.2) \qquad s_\eta = s + \tfrac{1}{2}\eta t s + \mathcal{o}(\eta)$$

holds in $L_2(\mu)$. Note that, in contrast to the definition in Bickel, Klaassen, Ritov and Wellner (1993), the dominating measure $\mu$ may depend on the particular path and that hence we do not have to assume that our model $\mathcal{P}$ is dominated. Taking squares in (4.2) and integrating with respect to $\mu$ we obtain $\int t\,dP = 0$, which we denote by $t \in L_2^0(P)$.

Let $\mathcal{C}_P$ be a collection of paths $\mathcal{P}_P$ in $\mathcal{P}$ through $P$. By the tangent set $\dot{\mathcal{P}}^0$ we denote the set of all tangents $t$ generated by paths in $\mathcal{C}_P$. The closed linear span $[\dot{\mathcal{P}}^0]$ of $\dot{\mathcal{P}}^0$ is called the tangent space in $\mathcal{P}$ at $P$ generated by the collection $\mathcal{C}_P$ of paths $\mathcal{P}_P$. This tangent space is denoted by $\dot{\mathcal{P}} \subset L_2^0(P)$.

Let $\mathcal{B}$ be a Banach space with norm $\|\cdot\|_{\mathcal{B}}$ and consider a map $\nu$ from $\mathcal{P}$ to $\mathcal{B}$. We shall call $\nu: \mathcal{P} \to \mathcal{B}$ pathwise differentiable at $P$ with respect to $\mathcal{C}_P$ if there exists a continuous, linear map $\dot{\nu}: \dot{\mathcal{P}} \to \mathcal{B}$ such that for every path $\mathcal{P}_P = \{P_\eta : |\eta| < 1\}$ in $\mathcal{C}_P$ passing through $P$ with tangent $t$,

$$\|\nu(P_\eta) - \nu(P) - \eta\dot{\nu}(t)\|_{\mathcal{B}} = o(\eta). \tag{4.3}$$

Following Section 2 of van der Vaart (1991) and Section 5.2 of Bickel, Klaassen, Ritov and Wellner (1993), we define the efficient influence functions $\dot{\nu}_{b^*}$ of $\nu$ as follows: for $b^*$ in the dual space $\mathcal{B}^*$ of $\mathcal{B}$ (the space of all bounded linear functions from $\mathcal{B}$ to $\mathbb{R}$), the map $b^* \circ \dot{\nu}: \dot{\mathcal{P}} \to \mathbb{R}$ is linear and bounded. Hence, by the Riesz representation theorem there exists a unique element $\dot{\nu}_{b^*} \in \dot{\mathcal{P}}$ such that for every $t \in \dot{\mathcal{P}}$,

$$b^* \circ \dot{\nu}(t) = \langle \dot{\nu}_{b^*}, t \rangle = E\dot{\nu}_{b^*} t.$$

Here $\langle \cdot, \cdot \rangle$ denotes the inner product in $L_2^0(P)$ and $E$ denotes expectation with respect to $P$. Note that this definition of efficient influence function depends on $\dot{\mathcal{P}}$ and hence on the choice of $\mathcal{C}_P$.

From now on we take $\mathcal{P}$ to be a semiparametric model $\mathcal{P} = \{P_{\theta,G} : \theta \in \Theta, G \in \mathcal{G}\}$, as in (1.1), with $\Theta \subset \mathbb{R}^k$ open, and $\mathcal{G} \subset \mathcal{H}$, where $\mathcal{H}$ is a Hilbert space. Fix $G \in \mathcal{G}$ and let $\mathcal{C}_G$ be a collection of paths in $\mathcal{G} \subset \mathcal{H}$ through $G$. By $\dot{\mathcal{G}}$ we denote the tangent space in $\mathcal{G}$ at $G$ generated by $\mathcal{C}_G$, that is, the closed linear span of tangents at $G$ along a path in $\mathcal{C}_G$. We focus on estimation of Banach-valued parameters of the form $\nu = \nu(P_{\theta,G}) = \tilde{\nu}(G)$, where $\tilde{\nu}: \mathcal{G} \to \mathcal{B}$ is pathwise differentiable; that is, there exists a bounded linear operator $\dot{\tilde{\nu}}: \dot{\mathcal{G}} \to \mathcal{B}$ such that for all paths $\{G_\eta : |\eta| < 1\} \in \mathcal{C}_G$ with tangent $\dot{G}$ [cf. (4.1)],

$$\|\tilde{\nu}(G_\eta) - \tilde{\nu}(G) - \eta\dot{\tilde{\nu}}(\dot{G})\|_{\mathcal{B}} = o(\eta). \tag{4.4}$$

Again, for every $b^* \in \mathcal{B}^*$, the map $b^* \circ \dot{\tilde{\nu}}: \dot{\mathcal{G}} \to \mathbb{R}$ is linear and continuous and hence there exists a unique $\dot{\tilde{\nu}}_{b^*} \in \dot{\mathcal{G}}$ such that for every $\dot{G} \in \dot{\mathcal{G}}$

$$b^* \circ \dot{\tilde{\nu}}(\dot{G}) = \langle \dot{\tilde{\nu}}_{b^*}, \dot{G} \rangle_{\mathcal{H}}. \tag{4.5}$$

The elements $\dot{\tilde{\nu}}_{b^*}$ for $b^* \in \mathcal{B}^*$ are called the gradients of $\tilde{\nu}$; they are similar to the efficient influence functions of $\nu$, described earlier.

If the parametric submodel $\mathcal{P}_1 = \mathcal{P}_1(G)$ of our semiparametric model $\mathcal{P}$ is regular in the sense of Definition 2.2, its tangent space $\dot{\mathcal{P}}_1$ is defined to be the closed linear span $[\dot{l}_1]$ of the $k$-vector of score functions $\dot{l}_1 = \dot{l}_1(\theta, G)$. This agrees with the definition of tangent spaces in arbitrary statistical models [cf. (4.2)] by several choices of a collection $\mathcal{C}_\theta$ of paths, for example, $\mathcal{C}_\theta = \{\{P_{\theta+\eta e_i, G} \in \mathcal{P}_1(G) : -1 < \eta < 1\} : i = 1, \ldots, k\}$ with $e_i$, $i = 1, \ldots, k$, unit vectors.

By keeping $\theta$ fixed and by varying $G$ we get another submodel $\mathcal{P}_2 = \mathcal{P}_2(\theta)$. Given a collection $\mathcal{C}_G$ of paths within $\mathcal{P}_2(\theta)$, the tangent space $\dot{\mathcal{P}}_2$ at $P_{\theta, G}$ is defined as the closed linear span in $L_2^0(P_{\theta, G})$ of all functions $\tau \in L_2^0(P_{\theta, G})$ such that

$$(4.6) \qquad s(\theta, G_\eta) = s(\theta, G) + \tfrac{1}{2}\eta \tau s(\theta, G) + o(\eta),$$

in $L_2(\mu)$, for some path $\{G_\eta : |\eta| < 1\} \in \mathcal{C}_G$. Note again that $\dot{\mathcal{P}}_2$ depends on the choice of $\mathcal{C}_G$. We assume that $\mathcal{C}_G$ is chosen in such a way that for every path $\{P_\eta = P_{\theta+\eta\zeta, G_\eta} : |\eta| < 1\}$ with $\{G_\eta : |\eta| < 1\} \in \mathcal{C}_G$, there exists a tangent $\rho \in L_2^0(P_{\theta, G})$ satisfying

$$(4.7) \qquad s(\theta + \eta\zeta, G_\eta) = s(\theta, G) + \tfrac{1}{2}\eta \rho s(\theta, G) + o(\eta),$$

in $L_2(\mu)$. The tangent space $\dot{\mathcal{P}}$ at $P_{\theta, G}$ is the closed linear span in $L_2^0(P_{\theta, G})$ of all these tangents $\rho \in L_2^0(P_{\theta, G})$. Typically, we have $\dot{\mathcal{P}} = [\dot{l}_1] + \dot{\mathcal{P}}_2$.

In fact, we will assume that the tangents from (4.7) have a special but frequently occurring structure, namely that of Hellinger differentiability.

DEFINITION 4.1. For every $\theta \in \Theta$ and $G \in \mathcal{G}$, the model $\mathcal{P}$ is *Hellinger differentiable* at $P_{\theta, G}$ if there exists a bounded linear operator $\dot{l} : \mathbb{R}^k \times \dot{\mathcal{G}} \to L_2^0(P_{\theta, G})$ such that for every $\zeta \in \mathbb{R}^k$ and every path $\{G_\eta : |\eta| < 1\} \in \mathcal{C}_G$ with tangent $\dot{G} \in \dot{\mathcal{G}}$,

$$(4.8) \qquad s(\theta + \eta\zeta, G_\eta) = s(\theta, G) + \tfrac{1}{2}\eta (\dot{l}(\zeta, \dot{G})) s(\theta, G) + o(\eta),$$

in $L_2(\mu)$.

The operator $\dot{l}$ is called the score operator. It may be expressed in terms of the score function $\dot{l}_1$ for $\theta$ in $\mathcal{P}_1(G)$ and the so-called score operator $\dot{l}_2$ for $G$ in $\mathcal{P}_2(\theta)$ as follows. For $\zeta \in \mathbb{R}^k$ and $\dot{G} \in \dot{\mathcal{G}}$, we have

$$(4.9) \qquad \dot{l}(\zeta, \dot{G}) = \dot{l}_1^\top \zeta + \dot{l}_2(\dot{G}).$$

Note that (4.8) and (4.9) reduce to (3.1) in the case where $\{G_\eta : |\eta| < 1\} = \{G\}$ is a singleton, and to (4.6) in case $\zeta = 0$.

In the following proposition we collect some fundamental results on the efficient influence functions for estimating $\theta$ in $\mathcal{P}$ and for estimating $\nu(P_{\theta, G}) =$ $\tilde{\nu}(G)$, both in the submodel $\mathcal{P}_2(\theta)$ and in the full model $\mathcal{P}$. The efficient influence function for estimating $\theta$ in $\mathcal{P}_1(G)$ is not of immediate interest for our purposes and hence is not discussed here. Define the efficient score function $l_1^*$ for estimating $\theta$ in the full model $\mathcal{P}$ by

$$(4.10) \qquad l_1^* = \dot{l}_1 - \Pi(\dot{l}_1 | \dot{\mathcal{P}}_2).$$

The efficient information matrix at $P_{\theta,G}$ for estimating $\theta$ in $\mathcal{P}$ is defined as

$$(4.11) \qquad I_*(\theta) = E(l_1^* l_1^{*T}).$$

Define the information operator as $\dot{l}_2^\top \dot{l}_2 : \dot{\mathcal{G}} \to \dot{\mathcal{G}}$ and let $(\dot{l}_2^\top \dot{l}_2)^- \alpha$ be a solution $h \in \dot{\mathcal{G}}$ of $\dot{l}_2^\top \dot{l}_2 h = \alpha$, for $\alpha \in \dot{\mathcal{G}}$. Let $N(A)$ and $R(A)$ denote the null space and the range of an operator $A$.

PROPOSITION 4.1. *Consider a map $\nu : \mathcal{P} \to \mathcal{B}$ given by $\nu(P_{\theta,G}) = \tilde{\nu}(G)$. Fix $\theta$, $G$ and $\mathcal{C}_G$, let $\mathcal{P}_1(\mathcal{G})$ be a regular parametric model as in Definition 2.2, and let $\mathcal{P}$ be Hellinger differentiable as in Definition 4.1. If:*

(i) *$\dot{\mathcal{G}}^0$ is a closed and linear subspace of $\mathcal{H}$, that is, $\dot{\mathcal{G}} = \dot{\mathcal{G}}_G = \dot{\mathcal{G}}^0$,*
(ii) *$\tilde{\nu} : \mathcal{G} \to \mathcal{B}$ is pathwise differentiable at $G$, as in (4.4),*
(iii) *$I_*(\theta)$ from (4.11) is nonsingular,*

*then*

A. *The efficient influence function at $P_{\theta,G}$ for estimating $\theta$ in $\mathcal{P}$ is given by*

$$(4.12) \qquad \tilde{l}_1 = I_*^{-1}(\theta) l_1^*.$$

B. *The map $\nu : \mathcal{P}_2(\theta) \to \mathcal{B}$ is pathwise differentiable at $P_{\theta,G}$ if and only if*

$$(4.13) \qquad \dot{\tilde{\nu}}_{b^*} \in R(\dot{l}_2^\top) \qquad \forall b^* \in \mathcal{B}^*.$$

*The efficient influence functions of $\nu$ are related to the gradients of $\tilde{\nu}$ by*

$$(4.14) \qquad \dot{\tilde{\nu}}_{b^*} = \dot{l}_2^\top \dot{\nu}_{b^*}, \qquad \dot{\nu}_{b^*} \in \dot{\mathcal{P}}_2, \dot{\tilde{\nu}}_{b^*} \in \dot{\mathcal{G}}.$$

*If also*

$$(4.15) \qquad \dot{\tilde{\nu}}_{b^*} \in R(\dot{l}_2^\top \dot{l}_2) \qquad \text{for all } b^* \in \mathcal{B}^*,$$

*then the unique solution of (4.14) is given by*

$$(4.16) \qquad \dot{\nu}_{b^*} = \dot{l}_2 (\dot{l}_2^\top \dot{l}_2)^- \dot{\tilde{\nu}}_{b^*}.$$

C. *The map $\nu:\mathcal{P} \to \mathcal{B}$ is pathwise differentiable at $P_{\theta,G}$ if and only if* (4.13) *holds. The efficient influence functions of $\nu$ are related to the gradients of $\tilde{\nu}$ by*

$$
\begin{aligned}
0 &= \langle \dot{l}_1, \dot{\nu}_{b^*} \rangle_\theta, \\
\dot{\tilde{\nu}}_{b^*} &= \dot{l}_2^\top \dot{\nu}_{b^*}.
\end{aligned}
\tag{4.17}
$$

*If also $\dot{\tilde{\nu}}_{b^*} \in R(\dot{l}_2^\top \dot{l}_2)$, for all $b^* \in \mathcal{B}^*$, then the unique solution of* (4.17) *is given by*

$$
\dot{\nu}_{b^*} = \dot{l}_2(\dot{l}_2^\top \dot{l}_2)^- \dot{\tilde{\nu}}_{b^*} - \langle \dot{l}_2(\dot{l}_2^\top \dot{l}_2)^- \dot{\tilde{\nu}}_{b^*}, \dot{l}_1 \rangle_\theta^\top I_*^{-1} l_1^*.
\tag{4.18}
$$

Parts B and C of this proposition are due to van der Vaart (1991); see his Theorem 3.1, formula (3.10) and Corollary 6.2. The gist of formula (4.18) is already contained in Begun, Hall, Huang and Wellner (1983), (4.4) and (3.1). Proofs of the proposition may be found also in Bickel, Klaassen, Ritov and Wellner (1993); see their Corollary 3.4.1, Theorem 5.4.1 and Corollaries 5.4.2 and 5.5.2. Note however, that they need the conditions $\dot{\mathcal{P}}_2 = \overline{R(\dot{l}_2)}$ and $\dot{\mathcal{P}} = \overline{R(\dot{l})}$. This is caused by their definition of tangent space $\dot{\mathcal{P}}$ as the closed linear span in $L_2^0(P_{\theta,G})$ of all possible tangents $\rho \in L_2^0(P_{\theta,G})$ that may be obtained via *some* path $\{P_\eta : |\eta| < 1, P_\eta \in \mathcal{P}\}$. In any particular model, the goal is construction of efficient estimators. The convolution theorem implies that if efficient estimators exist, they are asymptotically linear in the efficient influence functions; see Theorem 2.1 of van der Vaart (1991) and Theorems 3.3.2, 5.2.1 and 5.2.2 of Bickel, Klaassen, Ritov and Wellner (1993). In principle, the variances of the efficient influence functions corresponding to Bickel, Klaassen, Ritov and Wellner (1993) equal at least those corresponding to van der Vaart (1991), and should they differ, efficient estimators in the sense of van der Vaart (1991) do not exist. However, in practice estimators can be constructed that are efficient in this sense for appropriate choices of $\mathcal{C}$, which implies that they have to be efficient in the sense of Bickel, Klaassen, Ritov and Wellner (1993) as well. Of course, the advantage of the present approach is that the extra conditions mentioned above need not be verified now.

If also $N(\dot{l}_2) = \{0\}$ and $R(\dot{l}_2)$ is closed, then $\dot{l}_2^\top \dot{l}_2$ is one-to-one and onto, so $(\dot{l}_2^\top \dot{l}_2)^-$ may be replaced by $(\dot{l}_2^\top \dot{l}_2)^{-1}$. In this case all parameters $\nu(P)$ expressible as pathwise differentiable functions of $G$ are pathwise differentiable; see Corollary 3.3 of van der Vaart (1991).

**5. Efficient estimation of Banach parameters.** Let $X_1, \ldots, X_n$ be an i.i.d. sample from $P_{\theta,G} \in \mathcal{P}$, a semiparametric model as in (1.1). In this section we shall construct an efficient estimator of $\nu(P_{\theta,G}) = \tilde{\nu}(G) \in \mathcal{B}$ based on $X_1, \ldots, X_n$ within the model $\mathcal{P}$, using the constructions and the heredity of asymptotic linearity as studied in Section 2. As described in Section 1, we start with an efficient estimator of $\tilde{\nu}(G)$ within the submodel $\mathcal{P}_2(\theta)$, where $\theta$ is fixed and known and $G$ varies in $\mathcal{G}$. An estimator of $\nu(P_{\theta,G}) = \tilde{\nu}(G)$ within $\mathcal{P}_2(\theta)$ is of course allowed to depend on $\theta$. Let $\hat{\nu}_{\theta,n}$ be such a submodel estimator. In view of part B of Proposition 4.1 this estimator is efficient within the submodel $\mathcal{P}_2(\theta)$ with respect to the chosen collection $\mathcal{C}_G$ of paths if it is asymptotically linear in the efficient influence function given in (4.16) with the score operator $\dot{l}_2 = \dot{l}_2(\theta, G) : \dot{\mathcal{G}} \to L_2(P_{\theta,G})$ at $(\theta, G)$ depending on $\theta$ and $G$. Note that $\dot{l}_2$ depends on $\mathcal{C}_G$ since $\dot{\mathcal{G}}$ does. We shall need this asymptotic linearity locally uniformly in $\theta$, in the same way as in (3.3).

DEFINITION 5.1. Fix a subset $\mathcal{B}_0^*$ of $\mathcal{B}^*$, the dual space of the Banach space $\mathcal{B}$. The submodel estimator $\hat{\nu}_{\theta,n}$ is called $\mathcal{B}_0^*$-*weakly locally submodel efficient* at $P_{\theta_0,G}$ if for every sequence $\{\theta_n\}$ with $\{\sqrt{n}(\theta_n - \theta_0)\}$ bounded, and every $b^* \in \mathcal{B}_0^*$,

$$(5.1) \qquad \sqrt{n}\left|b^*(\hat{\nu}_{\theta_n,n} - \tilde{\nu}(G)) - \frac{1}{n}\sum_{i=1}^n \psi(X_i; \theta_n, G; b^*)\right| \overset{P_{\theta_0,G}}{\to} 0$$

holds with

$$(5.2) \qquad \psi(x; \theta, G; b^*) = [\dot{l}_2(\theta, G)(\dot{l}_2^\top(\theta, G)\dot{l}_2(\theta, G))^{-} \dot{\tilde{\nu}}_{b^*}](x).$$

The main result of our paper states that, under regularity conditions, if $\hat{\theta}_n = t_n(X_1, \ldots, X_n)$ is an efficient estimator of $\theta$ in $\mathcal{P}$ and if $\hat{\nu}_{\theta,n} = u_n(X_1, \ldots, X_n; \theta)$ is a weakly locally submodel efficient estimator of $\nu(P_{\theta,G}) = \tilde{\nu}(G)$ at $\theta_0$, then the substitution estimator $\hat{\nu}_{\hat{\theta}_n,n}$ is an efficient estimator of $\nu$ at $\theta_0$ in the semiparametric model $\mathcal{P}$; see the discussion in Section 4 after Proposition 4.1. If $\hat{\nu}_{\theta,n}$ is sufficiently smooth in $\theta$, this substitution estimator itself may be proved to be efficient; see Theorem 5.2 below. Without this extra condition we have to resort to a split-sample version of the substitution estimator, as in Section 2. Fix a sequence of integers $\{\lambda_n\}_{n=1}^\infty$ such that (2.7) holds, and define $\tilde{\theta}_{n1}$ and $\tilde{\theta}_{n2}$ as in (2.8). Analogously to (2.8) and (2.9), write

$$(5.3) \quad \tilde{\nu}^{(1)}_{\theta,\lambda_n} = u_{\lambda_n}(X_1, \ldots, X_{\lambda_n}; \theta), \qquad \tilde{\nu}^{(2)}_{\theta,n-\lambda_n} = u_{n-\lambda_n}(X_{\lambda_n+1}, \ldots, X_n; \theta)$$

and

$$(5.4) \qquad \hat{\nu}_n = \frac{\lambda_n}{n}\tilde{\nu}^{(1)}_{\tilde{\theta}_{n2},\lambda_n} + \frac{n-\lambda_n}{n}\tilde{\nu}^{(2)}_{\tilde{\theta}_{n1},n-\lambda_n}.$$

To prove efficiency of this estimator at $\theta_0 \in \Theta$ we will need the following smoothness condition, which is similar to (2.6). For every $h \in \dot{\mathcal{G}}$ and every sequence $\{\theta_n\}$ with $\{\sqrt{n}(\theta_n - \theta_0)\}$ bounded,

$$
\begin{aligned}
&\left| \frac{1}{\sqrt{n}} \sum_{i=1}^n \psi_h(X_i; \theta_n, G) \right. \\
&\left. - \frac{1}{\sqrt{n}} \sum_{i=1}^n \psi_h(X_i; \theta_0, G) - c_h(\theta_0)\sqrt{n}(\theta_n - \theta_0) \right| \overset{P_{\theta_0, G}}{\to} 0
\end{aligned}
\tag{5.5}
$$

and

$$c_h(\theta_n) \to c_h(\theta_0) \tag{5.6}$$

hold with

$$\psi_h(x; \theta, G) = \dot{l}_2(\theta, G)(\dot{l}_2^\top(\theta, G)\dot{l}_2(\theta, G))^- h(x) \tag{5.7}$$

and

$$c_h(\theta) = -E_\theta(\psi_h(X_1; \theta, G)\dot{l}_1^\top(\theta)(X_1)). \tag{5.8}$$

Furthermore, we write [cf. (5.2)]

$$c(\theta, G; b^*) = -E_\theta(\psi(X_1; \theta, G; b^*)\dot{l}_1^\top(\theta)(X_1)). \tag{5.9}$$

Lemmas 2.1 and 2.2 might be useful in checking conditions (5.5) and (5.6). Our main result is efficiency of $\hat{\nu}_n$ as follows.

THEOREM 5.1. *Fix $\theta_0 \in \Theta$ and $\mathcal{B}_0^* \subset \mathcal{B}^*$. Suppose that (5.5), (5.6) and the conditions of Proposition 4.1 are satisfied in model (1.1) for appropriately chosen collections $\mathcal{C}_G$ of paths. Suppose that the submodel estimator $\hat{\nu}_{\theta,n}$ is $\mathcal{B}_0^*$-weakly locally submodel efficient as in (5.1) and that (4.15) holds at $P_{\theta,G}$. If efficient estimation of $\theta$ is possible within $\mathcal{P}$ and if $\hat{\theta}_n$ is an efficient estimator of $\theta$ in $\mathcal{P}$, then $\hat{\nu}_n$ defined by (2.7), (2.8), (5.3) and (5.4) is a $\mathcal{B}_0^*$-weakly efficient estimator of $\nu$ from (1.2) within the full model $\mathcal{P}$ at $P_{\theta_0,G}$; that is, for every sequence $\{\theta_n\}$ with $\{\sqrt{n}(\theta_n - \theta_0)\}$ bounded and every $b^* \in \mathcal{B}_0^*$ [cf. (2.10), (4.18), (5.2) and (5.9)],*

$$
\begin{aligned}
&\sqrt{n}\left| b^*(\hat{\nu}_n - \nu(P_{\theta_n, G})) \right. \\
&\left. - \frac{1}{n} \sum_{i=1}^n [\psi(X_i; \theta_n, G; b^*) + c(\theta_n, G; b^*) I_*^{-1}(\theta_n) l_1^*(\theta_n)(X_i)] \right| \overset{P_{\theta_0, G}}{\to} 0.
\end{aligned}
\tag{5.10}
$$

PROOF. For every $b^* \in \mathcal{B}_0^*$, Theorem 2.1 may be applied and the local asymptotic linearity in (5.10) may be seen to yield efficiency via Proposition 4.1.C, (5.2) and (5.9). □

A closer look at the proofs of Theorems 5.1 and 2.1 with $\kappa'(\theta) = 0$ reveals that if the orthogonality

$$[\dot{l}_1(\theta_0)] \perp \dot{\mathcal{P}}_2(\theta_0) \tag{5.11}$$

holds, then $c(\theta_0, G; b^*)$ and the last term at the left-hand sides of (2.12) and (2.14) vanish, as does the second term at the right-hand side of (4.18). Hence it suffices for $\hat{\theta}_n$ to be $\sqrt{n}$-consistent at $\theta_0$ and we do not need (5.6) and (2.7), but instead

$$0 < \liminf_{n \to \infty} \frac{\lambda_n}{n} \leq \limsup_{n \to \infty} \frac{\lambda_n}{n} < 1. \tag{5.12}$$

We formulate this special case as a corollary.

COROLLARY 5.1 (Adaptive case). *Fix $\theta_0 \in \Theta$ and $\mathcal{B}_0^* \subset \mathcal{B}^*$. Suppose that the conditions of Proposition 4.1 are satisfied in model (1.1) for appropriately chosen collections $\mathcal{C}_G$ of paths and that for all sequences $\{\theta_n\}$ with $\{\sqrt{n}(\theta_n - \theta_0)\}$ bounded,*

$$\left| \frac{1}{\sqrt{n}} \sum_{i=1}^n \psi_h(X_i; \theta_n, G) - \frac{1}{\sqrt{n}} \sum_{i=1}^n \psi_h(X_i; \theta_0, G) \right| \overset{P_{\theta_0, G}}{\to} 0. \tag{5.13}$$

*Suppose furthermore that $\hat{\nu}_{\theta,n}$ is $\mathcal{B}_0^*$-weakly locally submodel efficient as in (5.1) and that (4.15) holds at $P_{\theta_0, G}$. If $\hat{\theta}_n$ is a $\sqrt{n}$-consistent estimator at $\theta_0$ and if the orthogonality (5.11) holds, then $\hat{\nu}_n$ defined by (5.4) and (5.12) is a weakly efficient estimator of $\nu$ from (1.2) within the full model $\mathcal{P}$ at $P_{\theta_0, G}$; that is, for every sequence $\{\theta_n\}$ with $\{\sqrt{n}(\theta_n - \theta_0)\}$ bounded and every $b^* \in \mathcal{B}_0^*$,*

$$\sqrt{n} \left| b^*(\hat{\nu}_n - \nu(P_{\theta_n, G})) - \frac{1}{n} \sum_{i=1}^n \psi(X_i; \theta_n, G; b^*) \right| \overset{P_{\theta_0, G}}{\to} 0. \tag{5.14}$$

REMARK 5.1. Our main result, Theorem 5.1, states that, assuming sufficient regularity of a semiparametric model $\mathcal{P}$, two conditions, namely efficiency of an estimator $\hat{\theta}_n$ of the finite-dimensional parameter in the full model $\mathcal{P}$, and submodel efficiency of an estimator $\hat{\nu}_{\theta,n}$ of a functional of the infinite-dimensional parameter $\nu$ in the submodel $\mathcal{P}_2(\theta)$ with $\theta$ fixed, are sufficient to guarantee efficiency of the combined estimator $\hat{\nu}_n = \hat{\nu}_{\hat{\theta}_n, n}$ in $\mathcal{P}$. The result derives from general expressions for the influence functions of substitution estimators of Section 2. These expressions can be used to pinpoint what is needed in terms of efficiency or what is allowed in terms of deviations from efficiency of the separate estimators $\hat{\theta}_n$ and $\hat{\nu}_{\theta,n}$ to achieve efficiency of the substitution estimator $\hat{\nu}_n$. Here we will derive conditions heuristically. Let $\hat{\theta}_n$ and $\hat{\nu}_{\theta,n}$ be asymptotically linear estimators with influence functions $\psi_1$ and $\psi_2$, respectively. Without loss of generality, they may be written as $\psi_1 = \tilde{l}_1 + \Delta_1$ and $\psi_2 = I_{22}^{-1}\dot{l}_2 + \Delta_2$, with $\Delta_1 \perp [\dot{l}_1, \dot{l}_2]$ and $\Delta_2 \perp \dot{l}_2$; see Proposition 3.3.1 of Bickel, Klaassen, Ritov and Wellner (1993). Then by Theorem 2.1 and (2.20) and (2.10) the influence function of the substitution estimator is given by

$$\psi_2 + c\psi_1 = I_{22}^{-1}\dot{l}_2 + \Delta_2 - (E((I_{22}^{-1}\dot{l}_2 + \Delta_2)\dot{l}_1^\top))(\tilde{l}_1 + \Delta_1)$$
$$= I_{22}^{-1}\dot{l}_2 - I_{22}^{-1}I_{21}\tilde{l}_1 + \Delta_2 - I_{22}^{-1}I_{21}\Delta_1 - (E(\Delta_2\dot{l}_1^\top))(\tilde{l}_1 + \Delta_1),$$

which equals the efficient influence function $\tilde{l}_2 = I_{22}^{-1}\dot{l}_2 - I_{22}^{-1}I_{21}\tilde{l}_1$ if and only if

$$(5.15) \qquad \Delta_2 = I_{22}^{-1}I_{21}\Delta_1 + (E(\Delta_2\dot{l}_1^\top))(\tilde{l}_1 + \Delta_1)$$

holds. If $\hat{\theta}_n$ is efficient, that is, if $\Delta_1 = 0$, (5.15) shows that we need $\Delta_2 = b^\top\tilde{l}_1$; so deviations from efficiency of $\hat{\nu}_{\theta,n}$ are permitted, provided they are in $[\tilde{l}_1]$, that is, provided these deviations are matrix multiples of $\hat{\theta}_n - \theta$. An example of this phenomenon is given in (6.26) in Example 6.4. If $\hat{\nu}_{\theta,n}$ is efficient, that is, if $\Delta_2$ vanishes, (5.15) reduces to $I_{22}^{-1}I_{21}\Delta_1 = 0$. This means that in the adaptive case ($I_{21} = 0$), $\hat{\theta}_n$ need not be efficient (see Corollary 5.1) and that in the nonadaptive case $\hat{\theta}_n$ has to be efficient in order to obtain efficiency of $\hat{\nu}_n$. Of course, also combinations of estimators are possible where neither of them is efficient, but in this case only a lucky shot might yield an efficient combined estimator $\hat{\nu}_n$.

REMARK 5.2. The first occurrences of the terminology "adaptive estimators" are in Beran (1974) and Stone (1975). In Pfanzagl and Wefelmeyer [(1982), pages 14 and 15], it is argued that this terminology is rather unfortunate since "adaptiveness" is a property of the model, namely (5.11) holds, and not of the estimators, which are just semiparametrically efficient. van Eeden (1970), who was the first to construct partially adaptive estimators of location in the one- and two-sample problem, calls her estimators efficiency-robust. Since the terminology of adaptiveness is quite common nowadays, we will stick to it, although Pfanzagl and Wefelmeyer (1982) are right, and we will call $\hat{\nu}_n$ of Corollary 5.1 an adaptive estimator of the Banach parameter $\nu(P)$.

REMARK 5.3. In the adaptive situation of the corollary the direct substitution estimator $\hat{\nu}_{\hat{\theta}_n,n}$ can also be shown to be efficient in the sense of (5.14) if $\hat{\theta}_n$ takes its values in a grid on $\mathbb{R}^k$ with meshwidth of the order $\mathcal{O}(n^{-1/2})$. This is the classical discretization technique of Le Cam (1956), which has also been used in our proof of Theorem 2.2.

The next theorem states that the direct substitution estimator $\hat{\nu}_{\hat{\theta}_n,n}$ is efficient in the general semiparametric model, if $\hat{\nu}_{\theta,n}$ is sufficiently smooth in $\theta$.

THEOREM 5.2. *Under the conditions of Theorem 5.1, let $\hat{\theta}_n$ be an efficient estimator of $\theta$; in the adaptive situation of (5.11) it suffices that $\hat{\theta}_n$ be $\sqrt{n}$-consistent. Fix $b^* \in \mathcal{B}_0^*$. If for all $\delta > 0$, $\varepsilon > 0$ and $c > 0$, there exist $\zeta > 0$ and $n_0 \in \mathbb{N}$ such that for all $n \geq n_0$,*

$$(5.16) \quad P_{\theta_0,G}\left(\sup_{\sqrt{n}|\theta-\theta_0|\leq c,\sqrt{n}|\theta-\tilde{\theta}|\leq \zeta} \sqrt{n}|b^*(\hat{\nu}_{\theta,n}-\hat{\nu}_{\tilde{\theta},n})| \geq \varepsilon\right) \leq \delta$$

*holds, then the substitution estimator $b^*\hat{\nu}_{\hat{\theta}_n,n}$ is an efficient estimator of $b^*\nu$ with $\nu$ from (1.2) within the full model $\mathcal{P}$ at $P_{\theta_0,G}$; that is, it satisfies (5.10).*

PROOF. Note that (5.16) is a translation of (2.15) and apply Theorem 2.2. □

REMARK 5.4. In the special case where $\mathcal{G}$ may be identified with a subset of Euclidean space, Theorem 5.1, Corollary 5.1 and Theorem 5.2 also apply. Here we give a heuristic argument why these results might be true in the Euclidean and hence the general case. Let $\Theta \subset \mathbb{R}^k$ and $\mathcal{H} \subset \mathbb{R}^l$ and let

$$\mathcal{P} = \{P_{\theta,\eta} : \theta \in \Theta, \eta \in \mathcal{H}\}$$

be a regular $(k+l)$-dimensional parametric model in the sense of Definition 2.2. We have identified $\mathcal{G}$ with $\mathcal{H}$ and hence we have $\dot{\mathcal{G}} = \mathbb{R}^l$, provided the class $\mathcal{C}$ of allowed paths is large enough. Define

$$\dot{l}_1(\theta,\eta) = \frac{\partial}{\partial \theta} \log p(x;\theta,\eta) \quad \text{and} \quad \dot{l}_2(\theta,\eta) = \frac{\partial}{\partial \eta} \log p(x;\theta,\eta)$$

as the score functions for $\theta$ and $\eta$, respectively, and the Fisher information matrix by

$$I(\theta,\eta) = \begin{pmatrix} I_{11}(\theta,\eta) & I_{12}(\theta,\eta) \\ I_{21}(\theta,\eta) & I_{22}(\theta,\eta) \end{pmatrix}$$

with $I_{ij}(\theta,\eta) = E\dot{l}_i \dot{l}_j^\top (\theta,\eta)$. Regularity of $\mathcal{P}$ implies that $I_{22}(\theta,\eta)$ and $I(\theta,\eta)$ are nonsingular. The efficient score function for estimating $\theta$ is given by

$$l_1^*(\theta,\eta) = \dot{l}_1(\theta,\eta) - I_{12}I_{22}^{-1}\dot{l}_2(\theta,\eta),$$

and with

$$I_*(\theta,\eta) = El_1^* l_1^{*\top}(\theta,\eta),$$

the efficient influence function for estimating $\theta$ is given by (cf. Proposition 4.1.A)

$$\tilde{l}_1(\theta,\eta) = I_*^{-1}(\theta,\eta)l_1^*(\theta,\eta).$$

We are interested in estimation of $\nu(\theta, \eta) = \tilde{\nu}(\eta)$ within $\mathcal{P}$. Let $\tilde{\nu} : \mathbb{R}^l \to \mathbb{R}^m$ be differentiable with $(m \times l)$ partial derivative matrix $\dot{\tilde{\nu}}$. Now $\dot{\tilde{\nu}}(\eta) I_{22}^{-1} \dot{l}_2(\theta, \eta)$ is the efficient influence function for estimating $\tilde{\nu}(\eta)$ in the submodel $\mathcal{P}_2(\theta)$. This coincides with formula (4.16) of Proposition 4.1.B. Note that the operator $\dot{l}_2 : \mathbb{R}^l \to L_2^0(P_{\theta,\eta})$ is represented by the column $l$-vector $\dot{l}_2$ via

$$\dot{l}_2(a) = a^\top \dot{l}_2, \qquad a \in \mathbb{R}^l, \tag{5.17}$$

that the operator $(\dot{l}_2^\top \dot{l}_2)^- : \mathbb{R}^l \to \mathbb{R}^l$ is represented by the nonsingular $(l \times l)$-matrix $I_{22}^{-1}(\theta, \eta)$, and that

$$\dot{\tilde{\nu}}_{b^*} = \dot{\tilde{\nu}}^\top(\eta) b^* \in \mathbb{R}^l, \qquad b^* \in \mathbb{R}^m. \tag{5.18}$$

According to formula (4.18) of Proposition 4.1.C the efficient influence function $\tilde{l}_\nu$ for estimating $\tilde{\nu}(\eta)$ in the full model $\mathcal{P}$ is given by

$$\tilde{l}_\nu(\theta, \eta) = \dot{\tilde{\nu}}(\eta) I_{22}^{-1} (\dot{l}_2(\theta, \eta) - I_{21} \tilde{l}_1(\theta, \eta)). \tag{5.19}$$

Fix $\theta_0$ and suppose that $\hat{\nu}_{\theta,n}$ is a (weakly) locally submodel efficient estimator of $\tilde{\nu}(\eta)$ within $\mathcal{P}_2(\theta_0)$, that is,

$$\sqrt{n} \left| \hat{\nu}_{\theta_n, n} - \nu - \frac{1}{n} \sum_{i=1}^n \dot{\tilde{\nu}}(\eta) I_{22}^{-1} \dot{l}_2(\theta_n, \eta)(X_i) \right| \overset{P_{\theta_0, \eta}}{\to} 0. \tag{5.20}$$

Substituting an estimator $\hat{\theta}_n$ of $\theta$ for $\theta_n$ with influence function $\psi$ under $\theta_0$ and using Taylor's expansion and the weak law of large numbers, we can formally argue as follows:

$$\begin{aligned}
&\sqrt{n}(\hat{\nu}_{\hat{\theta}_n, n} - \nu) \\
&\sim \frac{1}{\sqrt{n}} \sum_{i=1}^n \dot{\tilde{\nu}}(\eta) I_{22}^{-1} \dot{l}_2(\hat{\theta}_n, \eta)(X_i) \\
&\sim \frac{1}{\sqrt{n}} \sum_{i=1}^n \left[ \dot{\tilde{\nu}}(\eta) I_{22}^{-1} \dot{l}_2(\theta_0, \eta)(X_i) + \dot{\tilde{\nu}}(\eta) \frac{\partial}{\partial \theta} I_{22}^{-1} \dot{l}_2(\theta_0, \eta)(X_i)(\hat{\theta}_n - \theta_0) \right] \\
&\sim \frac{1}{\sqrt{n}} \sum_{i=1}^n \dot{\tilde{\nu}}(\eta) I_{22}^{-1} \dot{l}_2(\theta_0, \eta)(X_i) \\
&\quad + \frac{1}{n} \sum_{i=1}^n \dot{\tilde{\nu}}(\eta) \frac{\partial}{\partial \theta} I_{22}^{-1} \dot{l}_2(\theta_0, \eta)(X_i) \sqrt{n}(\hat{\theta}_n - \theta_0) \\
&\sim \frac{1}{\sqrt{n}} \sum_{i=1}^n \dot{\tilde{\nu}}(\eta) I_{22}^{-1} \dot{l}_2(\theta_0, \eta)(X_i) \\
&\quad + \dot{\tilde{\nu}}(\eta) E_{\theta_0} \frac{\partial}{\partial \theta} I_{22}^{-1} \dot{l}_2(\theta_0, \eta)(X_1) \frac{1}{\sqrt{n}} \sum_{i=1}^n \psi^\top(X_i).
\end{aligned}$$

By partial integration we have, under regularity conditions,

$$E_{\theta_0}\left(\frac{\partial}{\partial\theta}I_{22}^{-1}\dot{l}_2(\theta_0,\eta)(X_1)\right)$$
$$= \int\left(\frac{\partial}{\partial\theta}I_{22}^{-1}\dot{l}_2(\theta_0,\eta)\right)p(\theta_0,\eta)\,d\mu$$
$$= \frac{\partial}{\partial\theta}\int I_{22}^{-1}\dot{l}_2(\theta_0,\eta)p(\theta_0,\eta)\,d\mu - \int I_{22}^{-1}\dot{l}_2(\theta_0,\eta)\frac{\partial}{\partial\theta}p(\theta_0,\eta)\,d\mu$$
$$= -\int I_{22}^{-1}\dot{l}_2\dot{l}_1^\top p(\theta_0,\eta)\,d\mu = -I_{22}^{-1}I_{21}(\theta_0,\eta).$$

This means that the influence function of $\hat{\nu}_{\hat{\theta}_n,n}$ equals

(5.21) $$\dot{\tilde{\nu}}(\eta)I_{22}^{-1}(\dot{l}_2(\theta_0,\eta) - I_{21}(\theta_0,\eta)\psi),$$

which corresponds to $\tilde{l}_\nu(\theta_0,\eta)$ from (5.19) if $\psi = \tilde{l}_1$, that is, if $\hat{\theta}_n$ is efficient in $\mathcal{P}$.

The regularity of $\mathcal{P}$ and in particular continuity and nonsingularity of $I_{22}(\theta,\eta)$ imply (2.25), and hence Lemma 2.3 yields (5.5). Consequently, by Theorem 5.1 a split-sample modification of $\hat{\nu}_{\hat{\theta}_n,n}$ is an efficient estimator of $\nu$, if (5.20) is valid. By arguments as in Gong and Samaniego (1981) and under their extra regularity conditions it may be verified that the submodel maximum likelihood estimator $\hat{\nu}_{\theta,n}$ given $\theta$ satisfies both (5.20) and (5.16). Then Theorem 5.2 shows that $\hat{\nu}_{\hat{\theta}_n,n}$ is efficient if $\hat{\theta}_n$ is. Gong and Samaniego (1981) prove this directly and they call $\hat{\nu}_{\hat{\theta}_n,n}$ a pseudo maximum likelihood estimator.

**6. Examples.** In this section we shall present a number of examples that illustrate our main results, namely Theorem 5.1, Corollary 5.1 and Theorem 5.2. The first example expands on Example 3.1. The next examples are important semiparametric test-cases well known from textbooks; our results should in any case be applicable for those examples. Example 6.2 treats linear regression, which was used in Section 1 for motivation, for the particular case of a symmetric error distribution. For a possibly asymmetric error distribution we study the location problem in Example 6.4. These statistical models are parametrized linkage models, which are discussed in Example 6.3. A recurring theme in these examples is the idea that estimators based on residuals are actually estimators based on the unobservable errors, with the unknown parameter needed to construct these errors replaced by suitable estimators; see also Example 3.2. Example 6.5 considers the bootstrap, and we conclude in Example 6.6 with another well-known semiparametric model: the Cox proportional hazards model.

EXAMPLE 6.1 (Efficiency of sample variance). In Example 3.1 we have shown the local asymptotic linearity of the sample variance in the class of all distributions with finite fourth moment. At any point $P$ of this model $\mathcal{P}$ the tangent space is maximal and equals $\dot{\mathcal{P}} = L_2^0(P)$, provided the collection $\mathcal{C}_P$ of paths in $\mathcal{P}$ is chosen sufficiently large; see Example 3.2.1 of Bickel, Klaassen, Ritov and Wellner (1993) for an explicit construction, which is also valid in our more general framework. Consequently, any locally asymptotically linear estimator of the variance is efficient; see Theorem 3.3.1 of Bickel, Klaassen, Ritov and Wellner (1993). In particular, the sample variance is efficient. Of course, this conclusion can also be drawn from Theorem 5.2, since $n^{-1}\sum_{i=1}^n (X_i - \theta)^2$ is efficient within $\mathcal{P}_2(\theta)$ and $\bar{X}_n$ within $\mathcal{P}$ for the same reasons of linearity and maximal tangent spaces. This line of argument may be used to show efficiency of all sample central moments and, more generally still, for all functions $h$ with $n^{-1}\sum_{i=1}^n h(X_i - \bar{X}_n)$ estimating $\tilde{\nu}(G) = \int h\, dG$ within an appropriately broad class of distribution functions $G$.

EXAMPLE 6.2 (Symmetric error distribution in linear regression). Suppose we observe realizations of $X_i = (Y_i, Z_i)$, $i=1,\ldots,n$, which are i.i.d. copies of $X = (Y, Z)$. The random $k$-vector $Z$ and the random variable $Y$ are related by

$$Y = \theta^\top Z + \epsilon, \tag{6.1}$$

where $\epsilon$ is independent of $Z$ and symmetrically distributed about 0 with unknown distribution function $G$ and density $g$ with respect to Lebesgue measure $\lambda$. For deriving lower bounds we assume that $Z$ has known distribution $F$ and that $EZZ^\top$ is nonsingular. Note that the unknown Euclidean parameter $\theta \in \mathbb{R}^k$ is identifiable via

$$\theta = (EZZ^\top)^{-1} E(Zm(Y|Z)), \tag{6.2}$$

where $m(Y|Z)$ denotes the median of the conditional symmetric distribution of $Y$ given $Z$. We are interested in estimating the symmetric error distribution $\nu(P_{\theta,G}) = \tilde{\nu}(G) = G$.

The density of $X$ with respect to $\lambda \times F$ is given by

$$p(x; \theta, G) = p(y, z; \theta, G) = g(y - \theta^\top z). \tag{6.3}$$

We assume that $G$ has finite Fisher information $I(G) = \int (g'/g)^2 g\, d\lambda$ for location, and hence we have

$$\dot{l}_1(\theta)(X) = \dot{l}_1(X; \theta, G) = -Z\frac{g'}{g}(Y - \theta^\top Z) = -Z\frac{g'}{g}(\epsilon) \tag{6.4}$$

and

$$\mathcal{G} = \left\{ G \in L_\infty(\lambda) : g \geq 0, \int g\, d\lambda = 1,\ g(-\cdot) = g(\cdot),\ I(G) < \infty \right\}. \tag{6.5}$$

We embed $\mathcal{G}$ into $\mathcal{H} = L_2(\lambda)$ by taking square roots of densities. The Fisher information $I(\cdot)$ for location is lower semicontinuous on $\mathcal{G}$. Therefore, we will restrict $\mathcal{C}_G$ to those paths on which $I(\cdot)$ is continuous. Such paths may be constructed in the same way as at the end of Example 3.2.1 of Bickel, Klaassen, Ritov and Wellner (1993). Then we have, embedding $\dot{\mathcal{G}}$ into $L_2^0(G)$,

$$
\begin{aligned}
\dot{\mathcal{G}}^0 &= \{h \in L_2^0(G) : h(-\cdot) = h(\cdot), h' \in L_2^0(G)\}, \\
\dot{\mathcal{G}} &= \{h \in L_2^0(G) : h(-\cdot) = h(\cdot)\}.
\end{aligned}
\tag{6.6}
$$

Note that $\dot{l}_2(\theta)$ is the embedding of $\dot{\mathcal{G}}$ into $L_2^0(P_{\theta,G})$ given by

$$h \mapsto h(Y - \theta^\top Z), \tag{6.7}$$

whence $\dot{l}_2(\theta)(\dot{\mathcal{G}}) = \dot{\mathcal{P}}_2$. The finiteness and positivity of the Fisher information $I(G)$, the nonsingularity of $EZZ^\top$, the choice of $\mathcal{C}_G$, the $L_2$-continuity theorem for translations and (6.7) ensure regularity and Hellinger differentiability as described in Definitions 2.2 and 4.1, respectively. Furthermore, the symmetry of $h \in \dot{\mathcal{G}}$ and antisymmetry of $\dot{l}_1(\theta)$ imply

$$\dot{l}_1(\theta) \perp \dot{l}_2(\theta) h. \tag{6.8}$$

Thus we are in an adaptive situation here.

The map $\tilde{\nu} : \mathcal{G} \to \mathcal{B}$, the cadlag functions on $[-\infty, \infty]$ with sup-norm, is pathwise differentiable at $G \in \mathcal{G}$ with derivative [cf. Example 5.3.3, page 193, of Bickel, Klaassen, Ritov and Wellner (1993)]

$$\dot{\tilde{\nu}}(h)(t) = \int (\tfrac{1}{2}(\mathbf{1}_{(-\infty,t]}(x) + \mathbf{1}_{(-\infty,t]}(-x)) - G(t))h(x)\,dG(x). \tag{6.9}$$

Note that (6.6), (6.8) and (6.9) imply that the conditions of Proposition 4.1 are satisfied. Furthermore, the $L_2$-continuity theorem for translations implies (2.25). Consequently, Lemma 2.3 shows the validity of (5.13). Finally, note that (4.15) holds since $R(\dot{l}_2^\top(\theta)\dot{l}_2(\theta)) = \dot{\mathcal{G}}$ and $\dot{\tilde{\nu}}_{b^*} \in \dot{\mathcal{G}}$ by definition.

With $\theta$ known, an efficient estimator of $G$ is the symmetrized empirical distribution function of $\epsilon_1, \ldots, \epsilon_n$, given by

$$\hat{G}_{\theta,n}(x) = \tfrac{1}{2}(G_{\theta,n}(x) + \bar{G}_{\theta,n}(x)), \tag{6.10}$$

where

$$G_{\theta,n}(x) = n^{-1} \sum_{i=1}^n \mathbf{1}_{[\epsilon_i \leq x]} = n^{-1} \sum_{i=1}^n \mathbf{1}_{[Y_i - \theta^\top Z_i \leq x]} \tag{6.11}$$

and

$$\bar{G}_{\theta,n}(x) = 1 - \lim_{y \searrow x} G_{\theta,n}(-y) \tag{6.12}$$

[cf. Example 5.3.3, pages 193–195 of Bickel, Klaassen, Ritov and Wellner (1993)]. We note that $\hat{G}_{\theta,n}$ is weakly locally submodel efficient, since it is exactly linear in the efficient influence function; see just above (5.3.10), page 194 of Bickel, Klaassen, Ritov and Wellner (1993). Finally, by a method of Scholz (1971) we know that the maximum likelihood estimator of $\theta$ corresponding to the logistic density exists under any density within our model. Furthermore, this pseudo maximum likelihood estimator is $\sqrt{n}$-consistent [cf. Example 7.8.2, page 401, of Bickel, Klaassen, Ritov and Wellner (1993)]. In fact, efficient and hence adaptive estimators of the regression parameter $\theta$ have been constructed, for example, by Dionne (1981), Bickel (1982) and Koul and Susarla (1983).

Consequently, by Corollary 5.1 the split-sample estimator defined by (5.4), (6.10)–(6.12) and (5.12) is efficient. Note that this efficient estimator does not use any knowledge about the distribution of $Z$ and hence is also adaptive with respect to the distribution of $Z$.

Clearly, in practice one would not apply sample splitting, but use $\hat{G}_{\hat{\theta}_n,n}$ itself, which is the symmetrized empirical distribution function based on the residuals $\hat{\epsilon}_i = Y_i - \hat{\theta}_n^\top Z_i$. This yields an efficient estimator of $G$ if $\hat{\theta}_n$ is discretized as described in Remark 5.3. Without discretization $\hat{G}_{\hat{\theta}_n,n}$ is weakly efficient in the sense of Theorem 5.2 for most $b^* \in \mathcal{B}^*$, including the evaluation maps. To see this it suffices to verify (5.16) for empirical distributions of regression residuals, as is done in the following lemma.

LEMMA 6.1. *In the regression model* (6.1) *let the error have bounded density $g$ (not necessarily symmetric) and let $E|Z|$ be finite. Let $b^* \in \mathcal{B}^*$ be such that there exists a finite signed measure $\mu$ with $b^*(b) = \int b(x)\, d\mu(x)$ and $\|b^*\| = |\mu|([-\infty,\infty]) < \infty$. For such $b^*$, the smoothness condition* (5.16) *holds for*

$$(6.13) \qquad \hat{\nu}_{\theta,n} = \frac{1}{n}\sum_{i=1}^n \mathbf{1}_{[Y_i-\theta^\top Z_i,\infty)}(\cdot).$$

PROOF. Let $b^*$ be given and let $\mu$ be the corresponding signed measure with $\|b^*\| = C$. Let the density $g$ of $\epsilon_i$ be bounded by $B$ and assume $E|Z| = A$. For $\tilde{\eta} \in \mathbb{R}^k$ and $\lambda_n \to \infty$, $\lambda_n/\sqrt{n} \to 0$, Markov's inequality yields (note $e^z - 1 < 2z$ for $0 < z$ sufficiently small)

$$P\left(\int \frac{1}{\sqrt{n}}\sum_{i=1}^n \mathbf{1}_{[|\epsilon_i+\tilde{\eta}^\top Z_i - x| \leq 2\zeta|Z_i|/\sqrt{n}]}\, d|\mu|(x) \geq \varepsilon\right)$$

$$\leq E\exp\left(\lambda_n\left\{\frac{1}{\sqrt{n}}\sum_{i=1}^n \int \mathbf{1}_{[|\epsilon_i+\tilde{\eta}^\top Z_i - x| \leq 2\zeta|Z_i|/\sqrt{n}]}\, d|\mu|(x) - \varepsilon\right\}\right)$$

$$= \exp\Big\{n\log\Big(1 + E\Big(\exp\Big(\frac{\lambda_n}{\sqrt{n}}$$

(6.14)
$$\times \int \mathbf{1}_{[|\epsilon_1+\tilde{\eta}^\top Z_1-x|\leq 2\zeta|Z_1|/\sqrt{n}]}\, d|\mu|(x)\Big) - 1\Big)\Big)$$

$$- \varepsilon\lambda_n\Big\}$$

$$\leq \exp\Big\{2\lambda_n\sqrt{n}\int P\Big(|\epsilon_1+\tilde{\eta}^\top Z_1 - x| \leq \frac{2\zeta|Z_1|}{\sqrt{n}}\Big)\, d|\mu|(x) - \varepsilon\lambda_n\Big\}$$

$$\leq \exp\{(8ABC\zeta - \varepsilon)\lambda_n\} \qquad \text{as } n \to \infty.$$

Since $Y - \theta^\top Z \leq x < Y - \tilde{\theta}^\top Z$ implies $|Y - \theta^\top Z - x| \leq |\theta - \tilde{\theta}||Z|$, we obtain

(6.15)
$$\sup_{\sqrt{n}|\theta-\theta_0|\leq c, \sqrt{n}|\theta-\tilde{\theta}|\leq \zeta} \Big|\int \frac{1}{\sqrt{n}}\sum_{i=1}^n (\mathbf{1}_{[Y_i-\theta^\top Z_i \leq x]} - \mathbf{1}_{[Y_i-\tilde{\theta}^\top Z_i \leq x]})\, d|\mu|(x)\Big|$$

$$\leq \sup_{\sqrt{n}|\eta|\leq c} \int \frac{1}{\sqrt{n}}\sum_{i=1}^n \mathbf{1}_{[|\epsilon_i+\eta^\top Z_i-x|\leq \zeta|Z_i|/\sqrt{n}]}\, d|\mu|(x).$$

Consider the grid $\mathcal{G}_\zeta$ with meshwidth $2(kn)^{-1/2}\zeta$ with $k$ the dimension of $\theta$. By (6.15) and (6.14) the probability in (5.16) may be bounded by

(6.16)
$$P_{\theta_0,G}\Big(\sup_{\sqrt{n}|\tilde{\eta}|\leq c+\zeta, \tilde{\eta}\in\mathcal{G}_\zeta} \int \frac{1}{\sqrt{n}}\sum_{i=1}^n \mathbf{1}_{[|\epsilon_i+\tilde{\eta}^\top Z_i-x|\leq 2\zeta|Z_i|/\sqrt{n}]}\, d|\mu|(x) \geq \varepsilon\Big)$$

$$\leq \sum_{\sqrt{n}|\tilde{\eta}|\leq c+\zeta, \tilde{\eta}\in\mathcal{G}_\zeta} \exp\{(8ABC\zeta - \varepsilon)\lambda_n\}$$

$$\leq \Big(\frac{c}{\zeta}+1\Big)^k \exp\{(8ABC\zeta - \varepsilon)\lambda_n\},$$

which converges to 0 if $8ABC\zeta < \varepsilon$ holds. $\square$

We have proved the following result.

PROPOSITION 6.1. *Consider the linear regression model* (6.1) *with the covariate vector $Z$ and the error $\epsilon$ independent, both with unknown distributions. The matrix $EZZ^\top$ is nonsingular and the error distribution $G$ is assumed to be symmetric about zero with finite Fisher information for location. There exist $\sqrt{n}$-consistent and even adaptive estimators of $\theta$. For any such estimator, any estimator of $G$ defined by* (5.4), (6.10)–(6.12) *and* (5.12) *is weakly efficient in the sense of* (5.14), *that is, asymptotically linear in the efficient influence function given in* (6.9). *Furthermore, the direct substitution estimator $\hat{G}_{\hat{\theta}_n,n}$ is weakly efficient for all $b^* \in \mathcal{B}^*$ as in Lemma* 6.1.

REMARK 6.1. Note that with $k=1$ and $Z$ degenerate at 1, this proposition yields an efficient estimator of the error distribution in the classical symmetric location problem. The ordinary location problem will be treated in Example 6.4.

REMARK 6.2. The idea of using the residuals to assess the error distribution is quite standard and has been around for a long time, for instance in testing for normality.

REMARK 6.3. Interest might be in the standardized symmetric error distribution, that is, in $G$ standardized to have unit variance, as in Example 3.2. This leads to a nonadaptive situation in which approaches as in the next example should lead to efficient estimators.

EXAMPLE 6.3 (Parametrized linkage models). As in Example 3.2 we consider the statistical model of $n$ i.i.d. copies of a random variable $X$ that is linked to an error variable $\epsilon$ with distribution function $G$ via

$$t_\theta(X) = \epsilon,$$

with $t_\theta : \mathcal{X} \to \mathbb{R}$ measurable and $\theta \in \Theta \subset \mathbb{R}^k$. Let $\theta$ be given. The empirical distribution function of $t_\theta(X_i)$, $i=1,\ldots,n$, is (asymptotically) linear in the influence function

$$(6.17) \qquad x \mapsto \mathbf{1}_{[t_\theta(x) \leq \cdot]} - G(\cdot).$$

This influence function and hence the empirical distribution function itself are efficient in estimating the distribution function $G$ if $G$ and $\mathcal{G}$ are unrestricted.

Typically, however, $G$ is constrained to be symmetric (as in the preceding example) or to have, for example, mean 0. In general, if the constraints can be described by

$$\int \gamma \, dG = 0$$

for some fixed measurable function $\gamma : \mathbb{R} \to \mathbb{R}^l$, then the efficient influence function in estimating $G$ may be obtained from (6.17) by projection [cf., e.g., (6.2.6) of Bickel, Klaassen, Ritov and Wellner (1993)] and equals

$$(6.18) \quad x \mapsto \mathbf{1}_{[t_\theta(x) \leq \cdot]} - G(\cdot) - E(\mathbf{1}_{[\epsilon \leq \cdot]} \gamma^\top(\epsilon)) \{E\gamma(\epsilon)\gamma^\top(\epsilon)\}^{-1} \gamma(t_\theta(x)).$$

Under appropriate regularity conditions,

$$\hat{G}_{\theta,n}(t) = \frac{1}{n} \sum_{i=1}^{n} \mathbf{1}_{[t_\theta(X_i) \leq t]}$$

$$
(6.19) \quad -\left(\frac{1}{n}\sum_{i=1}^{n}\mathbf{1}_{[t_\theta(X_i)\leq t]}\gamma^\top(t_\theta(X_i))\right)
$$

$$
\times\left\{\frac{1}{n}\sum_{i=1}^{n}\gamma(t_\theta(X_i))\gamma^\top(t_\theta(X_i))\right\}^{-1}\frac{1}{n}\sum_{i=1}^{n}\gamma(t_\theta(X_i))
$$

is an efficient estimator of $G(t)$, $t \in \mathbb{R}$, within this restricted class $\mathcal{G}$ of constrained distribution functions, given $\theta$. Subsequently, a weakly efficient estimator of $G$ within the semiparametric model with $\theta$ unknown may be obtained via the theorems of Section 5.

We will present the details of this approach for the particular case of $k=1$, $Z=1$ a.s., that is, for the location model, in the next example.

EXAMPLE 6.4 (Error distribution in location problem). Let $X_1,\ldots,X_n$ be i.i.d. random variables, which are copies of a random variable $X$ with unknown distribution $P \in \mathcal{P}$ and distribution function $F$ on $\mathbb{R}$. It is well known that the empirical distribution function $\hat{F}_n$ is efficient in estimating $F$, when $F$ is completely unknown. Let us assume now that the $X_i$ have finite variance and mean $\theta$. It is well known also that the sample mean $\bar{X}_n$ is efficient in estimating $\theta$. With $t_\theta(X) = X - \theta = \epsilon$, the error distribution function $G \in \mathcal{G}$, the class of all distribution functions with mean zero, satisfies

$$G(t) = F(t+\theta), \qquad t \in \mathbb{R}.$$

Given $\hat{F}_n$ and $\bar{X}_n$, a natural estimator of the unknown error distribution $G$, which has mean zero, would be

$$(6.20) \qquad \hat{G}_n(t) = \hat{F}_n(t+\bar{X}_n), \qquad t \in \mathbb{R}.$$

In fact, $\hat{G}_n$ is an asymptotically efficient estimator of the error distribution function $G$, as may be shown by computation of the efficient influence function along the lines of Example 5.3.8 of Bickel, Klaassen, Ritov and Wellner (1993). Let $\Psi$ be a collection of bounded functions $\psi : \mathbb{R} \to \mathbb{R}$ with bounded uniformly continuous derivative $\psi'$ and let $\nu$ map $\mathcal{P}$ into the Banach space $l^\infty(\Psi)$ of bounded functions on $\Psi$ with the supremum norm such that

$$(6.21) \quad \nu(P_{\theta,G})(\psi) = \tilde{\nu}(G)(\psi) = G(\psi) = \int \psi(t)\,dG(t), \qquad \psi \in \Psi.$$

Thus, $G$ is identified via $\tilde{\nu}(G)$ provided the class $\Psi$ is rich enough. Indeed, $\hat{G}_n$ from (6.20) is efficient, that is,

$$\sqrt{n}\left(\hat{G}_n(\psi) - G(\psi) - \frac{1}{n}\sum_{i=1}^{n}\tilde{l}(X_i)(\psi)\right) = o_P(1), \qquad \psi \in \Psi,$$

holds with the efficient influence function $\tilde{l}$ equal to [cf. (6.18)]

$$\tilde{l}(x)(\psi) = \psi(x - \theta) - E_{P_{\theta,G}}\psi(X - \theta) - E_{P_{\theta,G}}\psi'(X - \theta)(x - \theta), \tag{6.22}$$
$$\psi \in \Psi.$$

If we apply the approach of Section 5, we need an efficient estimator of $G$ for the case where $\theta$ is known. As explained via the parametrized linkage models of Examples 3.2 and 6.3, the naive empirical distribution of $X_i - \theta$, $i = 1, \ldots, n$, is not efficient, but an explicit weighted empirical of the $X_i - \theta$, $i = 1, \ldots, n$, as given in the following proposition, is asymptotically efficient.

PROPOSITION 6.2 (Location known). *Let $X_1, \ldots, X_n$ be i.i.d. random variables with known mean $\theta$ and distribution function $G(\cdot - \theta)$, where $G$ is unknown with finite variance. The estimator*

$$\tilde{G}_{\theta,n}(t) = \frac{1}{n}\sum_{i=1}^{n}\left\{1 - \frac{(X_i - \theta)(\bar{X}_n - \theta)}{S_n^2(\theta)}\right\}\mathbf{1}_{(-\infty,t]}(X_i - \theta), \qquad t \in \mathbb{R}, \tag{6.23}$$

*with $S_n^2(\theta) = n^{-1}\sum_{i=1}^{n}(X_i - \theta)^2$, is weakly efficient in estimating the error distribution function $G$ with $G$ identified via $\nu: \mathcal{G} \to \mathcal{B} = l^{\infty}(L_2(G))$, $\nu(G)(\psi) = \int \psi \, dG$ for $\psi \in L_2(G)$.*

PROOF. Without loss of generality we may take $\theta = 0$. For $\psi$ square integrable with respect to $G$ we have

$$\tilde{G}_{0,n}(\psi) - G(\psi)$$
$$= \frac{1}{n}\sum_{i=1}^{n}\left\{\psi(X_i) - \int \psi \, dG - \frac{1}{nS_n^2(0)}\sum_{j=1}^{n}\psi(X_j)X_jX_i\right\} \tag{6.24}$$
$$= \frac{1}{n}\sum_{i=1}^{n}\left\{\psi(X_i) - \int \psi \, dG - \frac{\mathrm{cov}_G(\psi(X), X)}{\mathrm{var}_G X}X_i\right\} + o_P\left(\frac{1}{\sqrt{n}}\right),$$

where the last equality is implied by the law of large numbers. Consequently, $\tilde{G}_{0,n}$ is asymptotically linear in the efficient influence function as given in Example 6.2.1 of Bickel, Klaassen, Ritov and Wellner (1993); see (6.18). □

REMARK 6.4. Note that $\tilde{G}_{\theta,n}$ from (6.23) is just $\hat{G}_{\theta,n}$ from (6.19) for $t_\theta(x) = x - \theta$, written appropriately. $\tilde{G}_{\theta,n}$ is a signed measure; in its far tails it need not be monotone.

Plugging in $\theta = \bar{X}_n$ we obtain

$$\tilde{G}_{\bar{X}_n,n}(t) = \frac{1}{n}\sum_{i=1}^{n}\mathbf{1}_{(-\infty,t]}(X_i - \bar{X}_n) = \hat{F}_n(t + \bar{X}_n) = \hat{G}_n(t), \tag{6.25}$$

the estimator of $G$ from (6.20) which has been proved efficient above for $\mathcal{B} = l^\infty(\Psi)$ with a smaller set $\Psi$ than $L_2(G)$ as in Proposition 6.2. The sample splitting and substitution technique of Theorem 5.1 yields a different though similar efficient estimator of the error distribution $G$. Applying Lemma 6.1 and Theorem 5.2, we obtain the weak efficiency of $\hat{G}_n$ for another $\mathcal{B}$ and $\mathcal{B}^*$.

Note that plugging in $\bar{X}_n$ for $\theta$ into the empirical distribution function $\hat{F}_{\theta,n}(t) = \frac{1}{n}\sum_{i=1}^n \mathbf{1}_{(-\infty,t]}(X_i - \theta)$ of the $X_i - \theta$ yields the same estimator $\hat{G}_n(t)$ of $G(t)$. Although, as noted before, $\hat{F}_{\theta,n}(t)$ is not an efficient estimator for $\theta$ known, the combined estimator is. From Remark 5.1 we know that the substitution estimator $\hat{\nu}_{\hat{\theta}_n,n}$ can be efficient even if $\hat{\nu}_{\theta,n}$ is not efficient, as long as the influence function of $\hat{\nu}_{\theta,n}$ satisfies (5.15). In this case, this translates to $\hat{F}_{\theta,n}(t) = \tilde{G}_{\theta,n}(t) + b(\bar{X}_n - \theta) + \mathcal{O}_P(1)$, for every $t \in \mathbb{R}$ and some $b \in \mathbb{R}$. This is indeed the case, since by (6.23)

$$\begin{aligned}
\hat{F}_{\theta,n}(t) - \tilde{G}_{\theta,n}(t) &= \frac{\bar{X}_n - \theta}{S_n^2(\theta)} \cdot \frac{1}{n}\sum_{i=1}^n (X_i - \theta)\mathbf{1}_{(-\infty,t]}(X_i - \theta) \\
&= \frac{\bar{X}_n - \theta}{\sigma^2} \cdot (E(X - \theta)\mathbf{1}_{(-\infty,t]}(X - \theta) + \mathcal{O}_P(1)),
\end{aligned}$$
(6.26)

because of the law of large numbers.

The empirical likelihood approach of Owen (1991) has been applied by Qin and Lawless (1994) in their Example 3 (continued), page 314, to obtain another implicitly defined efficient estimator $G^*_{\theta,n}$ of $G$. $G^*_{\theta,n}$ is a proper distribution function and substitution of $\theta$ by $\bar{X}_n$ in $G^*_{\theta,n}$ yields $\hat{G}_n$ as well.

EXAMPLE 6.5 (Bootstrap). When constructing confidence intervals for the mean $\theta$ using the sample mean $\bar{X}_n = n^{-1}\sum_{i=1}^n X_i$, one needs the distribution of $\sqrt{n}(\bar{X}_n - \theta)$. It can be simulated once the distribution of $X - \theta = X_1 - \theta$ is known. By the fundamental rule of thumb of statistics this distribution of $X - \theta$ should be estimated when unknown. According to Example 6.4 an efficient estimator of this distribution is $\tilde{G}_{\bar{X}_n,n} = \hat{F}_n(\cdot + \bar{X}_n)$ from (6.25) and (6.20). In this way the distribution of $X - \theta$ under $F$ is estimated by the distribution of $X^*$, say, under $\hat{F}_n(\cdot + \bar{X}_n)$, which equals the distribution of $X^* - \bar{X}_n$ under $\hat{F}_n$. Via this approach we see why in the bootstrap world the distribution of $X - \theta$ under $F$ should be replaced by the distribution of $X^* - \bar{X}_n$ under $\hat{F}_n$.

EXAMPLE 6.6 (Baseline survival distribution in Cox's proportional hazards model). We observe i.i.d. copies of $X = (Z,T)$, where the hazard function of an individual with covariate $Z = z \in \mathbb{R}$ is given by

$$\lambda(t|z) = e^{\theta z}\lambda(t),$$

where $\theta \in \mathbb{R}$ and $\lambda$ is the so-called baseline hazard function, corresponding to covariate $z = 0$, and related to the Banach parameter $G$ as follows:

$$\lambda = \frac{g}{1-G} = \frac{g}{\overline{G}}. \tag{6.27}$$

Here $g$ is the density corresponding to the distribution function $G$ on $[0, \infty)$ of $T$, given $Z = 0$. Fix $T_0 > 0$ and define $\mathcal{G}$ to be all distribution functions $G$ with $G(T_0) < 1$. We assume that the distribution of $Z$ is known and has distribution function $F$. Furthermore, we denote Lebesgue measure on $(0, \infty)$ by $\mu$ and note that identification of $G$ with $\sqrt{g}$ yields $\mathcal{G} \subset L_2(\mu)$. Then the density of $(Z, T)$ with respect to $\mu \times F$ is

$$p(z, t; \theta, G) = e^{\theta z} g(t)(1 - G(t))^{(\exp(\theta z) - 1)}. \tag{6.28}$$

As in Example 3.4.2 of Bickel, Klaassen, Ritov and Wellner (1993), it is not difficult to see that

$$\dot{l}_1(z, t; \theta) = z(1 - e^{\theta z}\Lambda(t)) \tag{6.29}$$

with

$$\Lambda(t) = \int_0^t \lambda(s)\,ds = \int_0^t \frac{g(s)}{1-G(s)}\,ds = -\log(1-G(t)).$$

Representing $\dot{\mathcal{G}}$ in $L_2(G)$, we get $\dot{\mathcal{G}} = \dot{\mathcal{G}}^0 = L_2^0(G)$, and $\dot{l}_2 : \dot{\mathcal{G}} \to L_2(P_{\theta,g})$ is given by

$$(\dot{l}_2(\theta)a)(z, t) = a(t) + (e^{\theta z} - 1)\frac{\int_t^\infty a(s)\,dG(s)}{1 - G(t)}. \tag{6.30}$$

It is well known [cf. Tsiatis (1981)] that if $EZ^2 \exp(2\theta Z)$ is bounded uniformly in a neighborhood of $\theta_0$, then the Cox (1972) partial likelihood estimator $\hat{\theta}_n$ is (locally) regular and asymptotically linear in the efficient influence function $\|l_1^*\|^{-2} l_1^*$, where $l_1^*$ of (4.10) is given by

$$l_1^*(z, t; \theta) = \dot{l}_1(z, t; \theta) - \left(\frac{S_{1,\theta}}{S_{0,\theta}}(t) - e^{\theta z}\int_0^t \frac{S_{1,\theta}}{S_{0,\theta}}\,d\Lambda\right), \tag{6.31}$$

with

$$S_{i,\theta}(t) = E_\theta Z^i e^{\theta Z}\mathbf{1}_{[t,\infty)}(T), \qquad i = 0, 1. \tag{6.32}$$

A complete proof of efficiency in a strong sense is given in Klaassen (1989) under nondegeneracy and boundedness of $Z$.

We are interested in estimating the baseline distribution function $\tilde{\nu}(G) = G$ on an interval $[0, T_0]$ with $P_G(T > T_0) > 0$. In view of this bounded window we will restrict $\mathcal{C}_G$ to all paths at $G$ in $\mathcal{G}$ with tangent $h$ vanishing outside $[0, T_0]$, yielding

$$\dot{\mathcal{G}} = \{h \in L_2^0(G) : h = h\mathbf{1}_{[0,T_0]}\}. \tag{6.33}$$

Furthermore, we will assume that $|Z|$ is bounded a.s. by $C < \infty$. With the notation

$$\tilde{S}_{i,\theta}(t) = E_\theta Z^i e^{2\theta Z} \mathbf{1}_{[t,\infty)}(T), \qquad i = 0, 1, \tag{6.34}$$

we have

$$\dot{S}_{0,\theta}(t) = \frac{\partial}{\partial \theta} S_{0,\theta}(t) = S_{1,\theta}(t) - \tilde{S}_{1,\theta}(t)\Lambda(t). \tag{6.35}$$

To verify (5.5) we note that for $h \in \dot{\mathcal{G}}$ [cf. Example 6.7.1.A of Bickel, Klaassen, Ritov and Wellner (1993)]

$$\begin{aligned}
\psi_h(X;\theta) &= \dot{l}_2(\theta)(\dot{l}_2^\top(\theta)\dot{l}_2(\theta))^{-1} h(Z,T) \\
&= \left(\bar{G}h(T) - \int_T^\infty h\,dG\right)\Big/S_{0,\theta}(T) \\
&\quad - \int_0^T e^{\theta Z}\left(\bar{G}h(s) - \int_s^\infty h\,dG\right)\Big/S_{0,\theta}(s)\,d\Lambda(s)
\end{aligned} \tag{6.36}$$

holds. It follows from Example 3.5 in Schick (2001) that (2.25) holds for $\psi_\kappa = \psi_h$ from (6.36) where $h$ is associated with a $b^*$ corresponding to a signed measure $q$ on $[0, T_0]$ as in (6.40). Indeed, $h = \tilde{\nu}_{b^*}(x) = q([x, T_0]) - \int G\,d\mu$ holds, and

$$\begin{aligned}
&\bar{G}(t)h(t) - \int_t^\infty h(s)\,dG(s) \\
&\quad = \bar{G}(t)q([t,T_0]) - \int_t^\infty \mu([s,T_0])\,dG(s) = 0, \qquad t > T_0.
\end{aligned}$$

This yields (5.5) and (5.6) with $c_h(\theta) = E_\theta(\dot{l}_1(X,T;\theta)\psi_h(X,T;\theta))$.

Given the regression parameter $\theta$, the nonparametric maximum likelihood estimator $\hat{G}_{\theta,n}$ of the baseline distribution function $G$ may be derived from the nonparametric maximum likelihood estimator of the baseline cumulative hazard function $\Lambda$, as described in Section 1 of Johansen (1983), and it equals

$$\hat{G}_{\theta,n}(s) = 1 - \exp\left\{-\sum_{i=1}^n \mathbf{1}_{[0,s]}(T_i)\left(\sum_{j=1}^n \mathbf{1}_{[T_j \geq T_i]}e^{\theta Z_j}\right)^{-1}\right\}, \qquad s > 0. \tag{6.37}$$

Breslow [(1974), (7), page 93] proposed the Kaplan–Meier-type estimator

$$\tilde{G}_{\theta,n}(s) = 1 - \prod_{i=1}^n\left\{1 - \mathbf{1}_{[0,s]}(T_i)\left(\sum_{j=1}^n \mathbf{1}_{[T_j \geq T_i]}e^{\theta Z_j}\right)^{-1}\right\}. \tag{6.38}$$

Both these estimators are asymptotically linear in the efficient influence function

$$\begin{aligned}
\psi_{\mathbf{1}_{[0,s]}-G(s)}(z,t;\theta) &= \dot{l}_2(\theta)(\dot{l}_2^\top(\theta)\dot{l}_2(\theta))^{-1}(\mathbf{1}_{[0,s]}(\cdot) - G(s))(z,t) \\
&= \bar{G}(s)\left\{\frac{1}{S_{0,\theta}(t)}\mathbf{1}_{[0,s]}(t) - e^{\theta z}\int_0^{s\wedge t}\frac{1}{S_{0,\theta}}\,d\Lambda\right\},
\end{aligned} \tag{6.39}$$

uniformly in $s \in [0, T_0]$; see Section 4 of Tsiatis (1981) and Example 6.7.1.A of Bickel, Klaassen, Ritov and Wellner (1993). They are even weakly locally submodel efficient under the assumption of boundedness of $Z$ for $\mathcal{B}$ the cadlag functions on $[0, T_0]$ with supremum norm and $b^* \in \mathcal{B}_0^*$ of the type

$$b^*(b) = \int_{[0,T_0]} b(s) \, d\mu(s) \tag{6.40}$$

for some finite signed measure $\mu$. To verify this and for future use we need the following result.

LEMMA 6.2. *If $T_1, \ldots, T_n$ are random variables with empirical distribution function $\hat{F}_n$, then the statistic*

$$V_n(s) = \sum_{i=1}^n \mathbf{1}_{[T_i \leq s]} \left( \sum_{j=1}^n \mathbf{1}_{[T_j \geq T_i]} \right)^{-1} \tag{6.41}$$

*satisfies*

$$V_n(s) \leq -\log(1 - \hat{F}_n(s)). \tag{6.42}$$

PROOF. With $T_{(1)} \leq T_{(2)} \leq \cdots \leq T_{(n)}$ the order statistics we have

$$V_n(s) = \sum_{i=1}^n \mathbf{1}_{[T_{(i)} \leq s]} \frac{1}{n-i+1} = \sum_{i=1}^{n\hat{F}_n(s)} \frac{1}{n-i+1}$$
$$\leq \int_1^{n\hat{F}_n(s)+1} \frac{1}{n+1-x} \, dx = -\log(1 - \hat{F}_n(s)). \quad \square$$

We also need the following convergence result.

LEMMA 6.3. *Denote*

$$W_n(t; \theta) = \frac{1}{n} \sum_{j=1}^n \mathbf{1}_{[T_j \geq t]} Z_j e^{\theta Z_j} \left( \frac{1}{n} \sum_{j=1}^n \mathbf{1}_{[T_j \geq t]} e^{\theta Z_j} \right)^{-2}. \tag{6.43}$$

*In the Cox proportional hazards model of* (6.28) *with $|Z|$ bounded we have for $0 \leq s \leq T_0$,*

$$\frac{1}{n} \sum_{i=1}^n \mathbf{1}_{[T_i \leq s]} W_n(T_i; \theta) \overset{P_{\theta,G}}{\to} \int_0^s \frac{S_{1,\theta}}{S_{0,\theta}}(t) \, d\Lambda(t). \tag{6.44}$$

PROOF. Conditionally, given $T_i = t \leq s$, the statistic $W_n(t; \theta)$ converges in probability to $S_{1,\theta} S_{0,\theta}^{-2}(t)$ where both $W_n(t; \theta)$ and its limit are bounded a.s. Consequently, given $T_i \leq s$ the difference $|W_n(T_i; \theta) - S_{1,\theta} S_{0,\theta}^{-2}(T_i)|$ converges in mean to 0 and hence the lemma holds. $\square$

Combining these lemmata, we see that the nonparametric maximum likelihood estimator $\hat{G}_{\theta,n}(s)$ satisfies

$$\sqrt{n}(\hat{G}_{\theta_n,n}(s) - \hat{G}_{\theta,n}(s)) + \sqrt{n}(\theta_n - \theta)\bar{G}(s)\int_0^s \frac{S_{1,\theta}}{S_{0,\theta}}\,d\Lambda$$

$$= \sqrt{n}\int_\theta^{\theta_n}\left\{(\hat{G}_{\eta,n}(s) - 1)\frac{1}{n}\sum_{i=1}^n \mathbf{1}_{[T_i \leq s]}W_n(T_i;\eta)\right.$$

(6.45)
$$\left. + \bar{G}(s)\int_0^s \frac{S_{1,\theta}}{S_{0,\theta}}\,d\Lambda\right\}d\eta$$

$$= \mathcal{O}_P\left(\sqrt{n}\int_\theta^{\theta_n} \frac{1}{n}\sum_{i=1}^n \mathbf{1}_{[T_i \leq s]}|W_n(T_i;\eta) - W_n(T_i;\theta)|\,d\eta\right) + o_P(1)$$

$$= \mathcal{O}_P\left(\sqrt{n}\int_\theta^{\theta_n}(\eta - \theta)V_n(s)\,d\eta\right) + o_P(1) = o_P(1)$$

under $\theta$, uniformly in $s \in [0, T_0]$. Note that by (6.27), (6.30) and Lemma 2.2,

(6.46)
$$\frac{1}{\sqrt{n}}\sum_{i=1}^n \{\psi_{\mathbf{1}_{[0,s]} - G(s)}(Z_i, T_i; \theta_n) - \psi_{\mathbf{1}_{[0,s]} - G(s)}(Z_i, T_i; \theta)\}$$

$$+ \sqrt{n}(\theta_n - \theta)\bar{G}(s)\int_0^s \frac{S_{1,\theta}}{S_{0,\theta}}\,d\Lambda \overset{P_{\theta,G}}{\to} 0$$

holds uniformly in $s \in [0, T_0]$. The asymptotic linearity of $\hat{G}_{\theta,n}(s)$, (6.45) and (6.46) together imply that $\hat{G}_{\theta,n}(\cdot)$ is weakly locally submodel efficient on $[0, T_0]$ in the sense of (5.1) with $b^*$ as in (6.40). Finally, note [cf. (6.45)]

(6.47)
$$\sqrt{n}|\hat{G}_{\hat{\theta},n}(s) - \hat{G}_{\tilde{\theta},n}(s)| \leq \sqrt{n}\left|\int_{\tilde{\theta}}^{\hat{\theta}} \frac{1}{n}\sum_{i=1}^n \mathbf{1}_{[T_i \leq s]}W_n(T_i;\eta)\,d\eta\right|$$

$$\leq \sqrt{n}\left|\int_{\tilde{\theta}}^{\hat{\theta}} Ce^{3|\eta|C}\,d\eta\right|V_n(s).$$

By Lemma 6.2 this yields (5.16) with $b^* \in \mathcal{B}_0^*$ as in (6.40), since

(6.48)
$$P_\theta\left(\int_{[0,T_0]} V_n(s)d|\mu|(s) \geq c_0\varepsilon/\zeta\right)$$

$$\leq P_\theta(\hat{F}_n(T_0) \geq 1 - e^{-(c_0\varepsilon)/(|\mu|([0,T_0])\zeta)})$$

is arbitrarily small for $\zeta$ sufficiently small. We have proved that Theorem 5.2 may be applied and that the full nonparametric maximum likelihood estimator $\hat{G}_{\hat{\theta}_n,n}(s)$ of the baseline distribution function $G$ is efficient if $\hat{\theta}_n$ is efficient. By similar arguments this may be shown also for Breslow's estimator $\tilde{G}_{\hat{\theta}_n,n}(s)$.

PROPOSITION 6.3. *Consider the Cox proportional hazards model of* (6.28) *with the covariate $Z$ bounded a.s. in absolute value. If $\hat{\theta}_n$ is an efficient estimator of the regression parameter $\theta$, then both $\hat{G}_{\hat{\theta}_n,n}(s)$ and $\tilde{G}_{\hat{\theta}_n,n}(s)$ are weakly efficient in estimating $G(s)$ in the sense of* (5.10) *with $b^* \in \mathcal{B}_0^*$ as in* (6.40).

**Acknowledgment.** We are grateful for the many detailed suggestions of the referees, which led to significant improvements of the paper.

## REFERENCES

BEGUN, J. M., HALL, W. J., HUANG, W.-M. and WELLNER, J. A. (1983). Information and asymptotic efficiency in parametric-nonparametric models. *Ann. Statist.* **11** 432–452. MR696057

BERAN, R. (1974). Asymptotically efficient adaptive rank estimates in location models. *Ann. Statist.* **2** 63–74. MR345295

BICKEL, P. J. (1982). On adaptive estimation. *Ann. Statist.* **10** 647–671. MR663424

BICKEL, P. J., KLAASSEN, C. A. J., RITOV, Y. and WELLNER, J. A. (1993). *Efficient and Adaptive Estimation for Semiparametric Models*. Johns Hopkins Univ. Press, Baltimore. MR1245941

BOLTHAUSEN, E., PERKINS, E. and VAN DER VAART, A. W. (2002). *Lectures on Probability Theory and Statistics. Lecture Notes in Math.* **1781**. Springer, New York. MR1915443

BRESLOW, N. E. (1974). Covariance analysis of censored survival data. *Biometrics* **30** 89–99.

BUTLER, S. and LOUIS, T. (1992). Random effects models with non-parametric priors. *Statistics in Medicine* **11** 1981–2000.

COSSLETT, S. (1981). Maximum likelihood estimator for choice-based samples. *Econometrica* **49** 1289–1316. MR625785

COX, D. R. (1972). Regression models and life tables (with discussion). *J. Roy. Statist. Soc. Ser. B* **34** 187–220. MR341758

DIONNE, L. (1981). Efficient nonparametric estimators of parameters in the general linear hypothesis. *Ann. Statist.* **9** 457–460. MR606633

GILBERT, P. B., LELE, S. R. and VARDI, Y. (1999). Maximum likelihood estimation in semiparametric selection bias models with application to AIDS vaccine trials. *Biometrika* **86** 27–43. MR1688069

GONG, G. and SAMANIEGO, F. J. (1981). Pseudo maximum likelihood estimation: Theory and applications. *Ann. Statist.* **9** 861–869. MR619289

HUANG, J. (1996). Efficient estimation for the proportional hazards model with interval censoring. *Ann. Statist.* **24** 540–568. MR1394975

JOHANSEN, S. (1983). An extension of Cox's regression model. *Internat. Statist. Rev.* **51** 165–174. MR715533

KLAASSEN, C. A. J. (1987). Consistent estimation of the influence function of locally asymptotically linear estimates. *Ann. Statist.* **15** 1548–1562. MR913573

KLAASSEN, C. A. J. (1989). Efficient estimation in the Cox model for survival data. In *Proc. Fourth Prague Symposium on Asymptotic Statistics* (P. Mandl and M. Hušková, eds.) 313–319. Charles Univ., Prague. MR1051449

KLAASSEN, C. A. J. and PUTTER, H. (1997). Efficient estimation of the error distribution in a semiparametric linear model. In *International Symposium on Contemporary Multivariate Analysis and its Applications* (K. Fang and F. J. Hickernell, eds.) B.1–B.8. Hong Kong.

KLAASSEN, C. A. J. and PUTTER, H. (2000). Efficient estimation of Banach parameters in semiparametric models. Technical report, Korteweg-de Vries Institute for Mathematics.

KOUL, H. L. (1992). *Weighted Empiricals and Linear Models*. IMS, Hayward, CA. MR1218395

KOUL, H. L. (2002). *Weighted Empirical Processes in Dynamic Linear Models*, 2nd ed. Springer, New York. MR1911855

KOUL, H. L. and SUSARLA, V. (1983). Adaptive estimation in linear regression. *Statist. Decisions* **1** 379–400. MR736109

LE CAM, L. (1956). On the asymptotic theory of estimation and testing hypotheses. *Proc. Third Berkeley Symp. Math. Statist. Probab.* **1** 129–156. Univ. California Press, Berkeley. MR84918

LEHMANN, E. L. (1999). *Elements of Large-Sample Theory*. Springer, New York. MR1663158

LOYNES, R. M. (1980). The empirical distribution function of residuals from generalised regression. *Ann. Statist.* **8** 285–298. MR560730

MÜLLER, U. U., SCHICK, A. and WEFELMEYER, W. (2001). Plug-in estimators in semiparametric stochastic process models. In *Selected Proc. Symposium on Inference for Stochastic Processes* (I. V. Basawa, C. C. Heyde and R. L. Taylor, eds.) 213–234. IMS, Beachwood, OH. MR2002512

MURPHY, S. A., ROSSINI, A. J. and VAN DER VAART, A. W. (1997). Maximum likelihood estimation in the proportional odds model. *J. Amer. Statist. Assoc.* **92** 968–976. MR1482127

MURPHY, S. A. and VAN DER VAART, A. (2000). On profile likelihood (with discussion). *J. Amer. Statist. Assoc.* **95** 449–485. MR1803168

NIELSEN, G. G., GILL, R. D., ANDERSEN, P. K. and SØRENSEN, T. I. A. (1992). A counting process approach to maximum likelihood estimation in frailty models. *Scand. J. Statist.* **19** 25–43. MR1172965

OWEN, A. (1991). Empirical likelihood for linear models. *Ann. Statist.* **19** 1725–1747. MR1135146

PFANZAGL, J. and WEFELMEYER, W. (1982). *Contributions to a General Asymptotic Statistical Theory*. Springer, New York. MR675954

PIERCE, D. (1982). The asymptotic effect of substituting estimators for parameters in certain types of statistics. *Ann. Statist.* **10** 475–478. MR653522

QIN, J. and LAWLESS, J. (1994). Empirical likelihood and general estimating equations. *Ann. Statist.* **22** 300–325. MR1272085

RANDLES, R. (1982). On the asymptotic normality of statistics with estimated MR653521 parameters. *Ann. Statist.* **10** 462–474.

SCHICK, A. (2001). On asymptotic differentiability of averages. *Statist. Probab. Lett.* **51** 15–23. MR1820140

SCHOLZ, F. W. (1971). Comparison of optimal location estimators. Ph.D. dissertation, Univ. California, Berkeley.

STONE, C. (1975). Adaptive maximum likelihood estimators of a location parameter. *Ann. Statist.* **3** 267–284. MR362669

TSIATIS, A. A. (1981). A large sample study of Cox's regression model. *Ann. Statist.* **9** 93–108. MR600535

VAN DER VAART, A. W. (1991). On differentiable functionals. *Ann. Statist.* **19** 178–204. MR1091845

VAN EEDEN, C. (1970). Efficiency-robust estimation of location. *Ann. Math. Statist.* **41** 172–181. MR263194

Korteweg-de Vries Institute for Mathematics  
Universiteit van Amsterdam  
Plantage Muidergracht 24  
1018 TV Amsterdam  
The Netherlands  
e-mail: chrisk@science.uva.nl

Department of Medical Statistics  
Leiden University Medical Center  
P.O. Box 9604  
2300 RC Leiden  
The Netherlands  
e-mail: h.putter@lumc.nl